\documentclass[11pt,reqno]{amsart} 
\usepackage[margin=1in]{geometry}

\usepackage{algorithm}
\usepackage{algorithmicx}
\usepackage{algpseudocode}
\usepackage{verbatim}
\usepackage{blkarray}

\usepackage[english]{babel}
\usepackage{amsmath}
\usepackage{amssymb}
\usepackage{amsthm}
\usepackage{latexsym}
\usepackage{tikz-cd}
\usepackage{mathtools}
\usepackage{amsfonts}
\usepackage{mathrsfs}
\usepackage{textcomp}
\usepackage[T1]{fontenc}
\usepackage{graphicx}
\usepackage{setspace}
\usepackage{nicefrac}
\usepackage{indentfirst}
\usepackage{dcolumn}
\newcolumntype{d}[1]{D{.}{.}{#1}}
\usepackage{enumerate}
\usepackage{wasysym}
\usepackage{upgreek}
\usepackage{makecell}
\usepackage{hyperref}
\usepackage{booktabs}
\usepackage{etoolbox}
\usepackage{wrapfig}
\usepackage{floatflt}
\usepackage{parskip}
\usepackage{lscape}
\usetikzlibrary{automata}
\newtheorem{theorem}{Theorem}[section]
\newtheorem{lemma}[theorem]{Lemma}
\newtheorem{corollary}[theorem]{Corollary}
\newtheorem{proposition}[theorem]{Proposition}
\theoremstyle{definition}
\newtheorem{definition}[theorem]{Definition}
\newtheorem{remark}[theorem]{Remark}
\newcommand{\mydef}[1]{{\color{blue} #1}}

\usepackage{graphicx}

\renewcommand{\dim}{\text{dim}}

\newcommand{\p}{{\bf p}}

\newcommand{\x}{{\bf x}}

\newtheorem{example}[theorem]{Example}

\newcommand{\ZZZ}{{\mathbb{Z}/2\mathbb{Z}}}
\newcommand{\Klein}{{V}}
\newcommand{\DDD}{{D_8}}

\usetikzlibrary{positioning,arrows}
\usepackage{blkarray}
\newcommand{\alert}[1]{{\color{red}#1}}

\newcommand{\CHANGE}[2]{#2}

\makeatletter
\newcommand{\ostar}{\mathbin{\mathpalette\make@circled\star}}
\newcommand{\make@circled}[2]{%
  \ooalign{$\m@th#1\smallbigcirc{#1}$\cr\hidewidth$\m@th#1#2$\hidewidth\cr}%
}
\newcommand{\smallbigcirc}[1]{%
  \vcenter{\hbox{\scalebox{0.77778}{$\m@th#1\bigcirc$}}}%
}
\makeatother

\newcommand{\OrbitsX}[4]{ 
\begin{tikzpicture}[scale = 0.2]
\coordinate (v1) at (0,0); 
\coordinate (v2) at (6,0);
\coordinate (v4) at (6,6);
\coordinate (v3) at (0,6);

\draw[thick,gray!50] (v1)--(v3) (v2)--(v4)  (v1) -- (v2) -- (v3) -- (v4) -- cycle ;

\foreach \va/\vb in {#1}{
\draw[line width = 3.5,opacity=1.0, blue!75] (v\va) -- (v\vb); 
}

\foreach \va/\vb/\vc in {#2}{
\draw[ thick,dashed, draw=black, fill=black!10!red, opacity=0.23] (v\va) -- (v\vb)-- (v\vc) --cycle; 
}

\foreach \tet in {#3}{
\draw[fill=black!30!green] (3,3) circle (20pt);
}

\foreach \vt in {1,...,4}{
\node[draw, circle, inner sep=1.5pt, fill=black] at (v\vt) {}; 
}
\node[above ] at (3,6) {\large #4};
\end{tikzpicture}
}

\newcommand{\OrbitsXNonbasis}[4]{ 
\begin{tikzpicture}[scale = 0.2]
\coordinate (v1) at (0,0); 
\coordinate (v2) at (6,0);
\coordinate (v4) at (6,6);
\coordinate (v3) at (0,6);

\draw[thick,gray!50] (v1)--(v3) (v2)--(v4)  (v1) -- (v2) -- (v3) -- (v4) -- cycle ;

\foreach \va/\vb in {#1}{
\draw[line width = 3.5,opacity=1.0, gray!50] (v\va) -- (v\vb); 
}

\foreach \va/\vb/\vc in {#2}{
\draw[ thick,dashed, draw=black, fill=gray!50, opacity=0.40] (v\va) -- (v\vb)-- (v\vc) --cycle; 
}

\foreach \tet in {#3}{
\draw[fill=gray!50] (3,3) circle (20pt);
}

\foreach \vt in {1,...,4}{
\node[draw, circle, inner sep=1.5pt, fill=gray!50] at (v\vt) {}; 
}
\node[above ] at (3,6) {\large #4};
\end{tikzpicture}
}

\title{The Algebraic Matroid of the Heron Variety}

\author{Seth K. Asante}
\address{Theoretisch-Physikalisches Institut\\ Friedrich-Schiller-Universit\"at Jena\\ Max-Wien-Platz 1, 00743 Jena, Germany}
\email{seth.asante@uni-jena.de}

\author{Taylor Brysiewicz}
\address{Department of Mathematics\\ University of Western Ontario\\
 2004 Perth Dr, London, ON N6G 2V4, Canada}
\email{tbrysiew@uwo.ca}

\author{Michelle Hatzel}
\address{Department of Mathematics\\ University of Western Ontario\\
 2004 Perth Dr, London, ON N6G 2V4, Canada}
\email{mhatzel@uwo.ca}

\makeatother

\begin{document}

\begin{abstract}
We introduce the $n$-th Heron variety as the realization space of the (squared) volumes of faces of an $n$-simplex. Our primary goal is to understand the extent to which Heron's formula, which expresses the area of a triangle as a function of its three edge lengths, can be generalized. Such a formula for one face volume of an $n$-simplex in terms of other face volumes expresses a dependence in the algebraic matroid of the Heron variety. Whether the volume is expressible in terms of radicals is controlled by the \CHANGE{Galois}{monodromy} groups of the coordinate projections of the Heron variety onto coordinates of bases. We discuss a suite of algorithms, some new, for determining these matroids and  \CHANGE{Galois}{monodromy} groups. We apply these algorithms toward the smaller Heron varieties, organize our findings, and interpret the results in the context of our original motivation. 
\end{abstract} 

\maketitle

\vspace{-0.4in}
\section{Introduction}
\noindent Heron's formula expresses the area $v_{123}$ of a triangle $\Delta (p_1,p_2,p_3)$ via its edge lengths~$v_{ij}~=~||\overline{p_ip_j}||$:
\begin{equation}
\label{eq:Heron}
v_{123} = \frac{1}{4}\sqrt{2v_{12}^2v_{13}^2+2v_{12}^2v_{23}^2+2v_{13}^2v_{23}^2-v_{12}^4-v_{13}^4-v_{23}^4}.
\end{equation}
The volume $v_{12\cdots n+1}$ of an $n$-simplex $\Delta_n(p_1\cdots p_{n+1})$ is similarly expressed in terms of its edge~lengths:
\begin{equation}
\label{eq:CMformula}
v_{12\cdots n+1} = \sqrt{\frac{(-1)^{n+1}}{(n!)^2\cdot 2^{n}} \cdot  \begin{vmatrix}
0 & 1 & 1 & 1 & \cdots & 1 \\
1 & 0 & v_{1,2}^2 & v_{1,3}^2 & \cdots & v_{1,n+1}^2 \\
1 & v_{1,2}^2 & 0 & v_{2,3}^2 & \cdots & v_{2,n+1}^2 \\
\vdots & \vdots & \vdots & \vdots & \ddots & \vdots \\
1 & v_{1,n+1}^2 & v_{2,n+1}^2 & v_{3,n+1}^2  & \cdots & 0
\end{vmatrix}}.
\end{equation} The faces of $\Delta_n$ are themselves simplices and so their volumes, up to a square root and a scaling, are the principal minors of the matrix in \eqref{eq:CMformula}. Consequently, the realization space which encodes the incidence of volumes of all $2^{n+1}-(n+1)-1=:N(n)$  positive-dimensional faces of $\Delta_n$ is parametrized by the $e(n):= {{n+1}\choose{2}}$ edge lengths. Under the change of coordinates $x_{12\ldots n+1} = v_{12\ldots n+1}^2$, both Heron's formula and \CHANGE{the Cayley-Menger formula}{\eqref{eq:CMformula}} become polynomial equations in the squared edge lengths $x_{ij}=v_{ij}^2$. We define the $n$-th Heron variety $X_n$ to be the algebraic closure in $\mathbb{C}_{\textbf{x}}^{N(n)}$ of this realization space in the squared volume coordinates $\textbf{x}$.

The above discussion may be summarized by the following statement: \textit{The edge lengths of a simplex determine all of its face volumes.} We address a more general question:
\begin{center}
\textbf{Q1 (Identifiability):} \emph{Which sets of volumes of faces of }$\Delta_n$\emph{  determine all face volumes of }$\Delta_n$\emph{?}
\end{center}

We tackle this question from the point-of-view of \textit{algebraic matroids}. The \textit{bases} of the \textit{algebraic matroid} associated to an irreducible affine variety $X\subseteq \mathbb{C}^N$ are those subsets $S$ of coordinates for which the projection $\pi_S:X \to \mathbb{C}_{S}^{|S|}$ is finite-to-one and dominant. In other words, these are the subsets of coordinates for which a generic value $\textbf{x}_S \in \mathbb{C}_{S}^{|S|}$ may be completed to finitely many points on $X$. Equation \eqref{eq:CMformula} implies that the squared edge length coordinates $\{x_{ij}\}_{1 \leq i < j < n+1}$ form a basis of the algebraic matroid $\mathcal M_n$ of $X_n$. A slight relaxation of \textbf{Q1} becomes
\begin{center}
\textbf{Q2 (The Algebraic Matroid):} \emph{What are the bases of the algebraic matroid $\mathcal M_n$ of $X_n$?}
\end{center}

Given an irreducible variety $X \subseteq \mathbb{C}^N$ and a basis $B$ of its algebraic matroid, the projection map $\pi_B:X \to \mathbb{C}_{B}^{|B|}$ can be thought of as a branched cover. From this point-of-view, $\pi_B$ has a \textit{degree} $d_B$ and a \textit{monodromy group} $G_{\pi_B}$. Its degree $d_B$ is the number of points in a generic fibre and its monodromy group is  comprised of all permutations obtainable by analytically continuing the points in that fibre over a loop in $\mathbb{C}_{B}^{|B|}$. The group $G_{\pi_B} \subseteq \mathfrak S_{d_B}$ describes the symmetries of the fibres of $\pi_B$. In \cite{Harris} Harris showed that this monodromy group is the Galois group associated to the field extension of function fields  induced by $\pi_B$. Hence, the solvability of this group controls whether the unknown volumes are solvable by radicals in the known ones\CHANGE{}{.}
\begin{center}
\textbf{Q3 (Monodromy Groups):} \emph{For each basis $B$ of $\mathcal M_n$, what is the monodromy group $G_{\pi_B}$?}
\end{center}
The degree $d_B$ of the map $\pi_B$ is implicit in an answer to question $\textbf{Q3}$. The numbers $d_B$ associated to the bases of the algebraic matroid of a variety are called \textit{base degrees} and are considered \textit{decorations} of the matroid (see \cite{Rosen}). As far as we know, the more refined decoration of a basis by the monodromy/Galois group is a new consideration in this work. 

We address the questions  above from an experimental viewpoint by combining computational methods in group theory, polyhedral geometry, symbolic algebraic geometry, and numerical algebraic geometry. We fully describe the algebraic matroids $\mathcal M_2,\mathcal M_3$ and $\mathcal M_4$ and provide a deeper analysis of the Galois groups  associated to $\mathcal M_3$.  The \texttt{julia} \cite{julia} software packages \texttt{OSCAR.jl} \cite{OSCAR} and  \texttt{HomotopyContinuation.jl} \cite{HCjl} are crucial to our work. Additionally, we use the new package  \texttt{Pandora.jl} \cite{Pandora} which combines the power of both \texttt{OSCAR.jl} and \texttt{HomotopyContinuation.jl} to perform computations on  enumerative problems. Code to reproduce our results can be found at  \vspace{-5pt}
\begin{center}
{\footnotesize{{\color{blue}{\href{https://github.com/tbrysiewicz/CodeForPapers/tree/main/TheAlgebraicMatroidOfTheHeronVariety}{https://github.com/tbrysiewicz/CodeForPapers/tree/main/TheAlgebraicMatroidOfTheHeronVariety}}}}}
\end{center}
 \vspace{-5pt}
We give an impression of our main results through three examples involving the tetrahedron $\Delta_3$\CHANGE{:}{.}

\begin{example}
Consider three $e(3)=6$-subsets of the $N(3) = 11$ positive-dimensional faces of $\Delta_3$\CHANGE{.}{:} 
\[
 B_{9}=\{12,13,24,34,123,124\}\quad  B_{10} = \{12,13,24,34,123,234\}  \quad  B_{31}=\CHANGE{\{12,23,34,124,134,1234\}}{\{12,24,34,123,134,1234\}}
\]
\begin{center}
\includegraphics[scale=0.5]{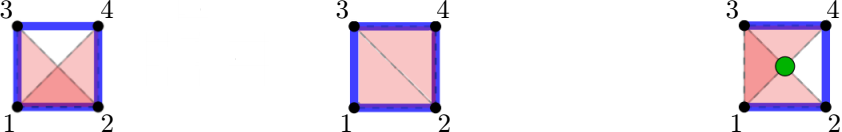}
\end{center}
The subsets $B_{9}$ and $B_{31}$ are bases of $\mathcal M_3$, but $B_{10}$ is not a basis. 
The monodromy groups of $\pi_{B_9}$ and $\pi_{B_{31}}$ are 
\[G_{\pi_{B_9}} = \ZZZ \times \ZZZ =: V \subseteq \mathfrak S_4 \quad \quad \text{ and } \quad \quad G_{\pi_{B_{31}}} = \mathfrak S_{12},
\]
respectively.
The group $V$ is solvable, but $\mathfrak S_{12}$ is not. As a consequence, there exist formulae in radicals for the volumes  $v_{23},v_{14},v_{134},v_{234},v_{1234}$ in terms of the volumes indexed by $B_9$, e.g.
\[v_{234} ={\scriptsize{ \frac{\sqrt{\left(-4 {v_{12}}^{2}-4 {v_{13}}^{2}+4 {v_{24}}^{2}+4 {v_{34}}^{2}\right) \sqrt{{v_{12}}^{2} {v_{13}}^{2}-4 {v_{123}}^{2}}-{v_{12}}^{4}+\cdots -{v_{34}}^{4}+16 {v_{123}}^{2}}}{4}}}.
\]
On the other hand, no such formula exists for any of the volumes \CHANGE{$v_{13},v_{14},v_{24},v_{123},v_{234}$}{$v_{13},v_{14},v_{23},v_{124},v_{234}$} in terms of those indexed by $B_{31}$. 
Similar statements for the remaining $32$-many $\mathfrak S_{4}$-orbits $B_1,\ldots,B_{35}$ (as displayed in Figure \ref{fig:X3Picture}) of $6$-subsets of positive-dimensional faces of $\Delta_3$ can be found in Section \ref{sec:results}.
\end{example}

\noindent \textbf{Relationship to previous and concurrent work:} 
The dependencies between the minors (principal or otherwise) of a matrix (symmetric, antisymmetric, or otherwise) of some rank (full or otherwise) is a popular type of question (e.g. \cite{Ahmadieh, Bruns, LinSturmfels, Oeding}). Our work fits into this framework, as we seek to understand the variables which support dependencies among the principal minors of a generic matrix of the form of the one in \eqref{eq:CMformula}.

A variety of edge-lengths $v_{ij}$ for which the matrix in \eqref{eq:CMformula} has some bounded rank is called a \textit{Cayley-Menger variety}. The algebraic matroids of the Cayley-Menger varieties, called \textit{rigidity matroids},  play a fundamental role in rigidity theory \cite{DanielBernstein1,Borcea2,Borcea1,Graver,Whiteley}. We leave the connection between the rigidity matroids and the algebraic matroids of the Heron varieties to future research.

The motivation for this project came from the question of whether or not the $10$-subset of triangular faces of $\Delta_4$ forms a basis for the rank $10$ matroid $\mathcal M_4$. The answer to this question has consequences in the world of theoretical physics where area volumes are used to (locally) identify $4$-simplices in a triangulation of space-time. This is the example in Section \ref{secsec:physics}, which is addressed in much greater detail in concurrent work by the first two authors \cite{SethTaylor}.

\noindent \textbf{Outline of paper:} Section \ref{sec:background} establishes the notation and language used throughout the paper. In Section \ref{sec:pipeline} we lay out the computational pipeline for addressing the questions $\textbf{Q1}-\textbf{Q3}$. We collect our computational results in Section \ref{sec:results} for $n<5$, using the algorithms of Section \ref{sec:pipeline}. In that same section, we also include some partial results for $n \geq 5$. Finally, in Section \ref{sec:experiments} we turn our attention to semi-algebraic questions about the Heron variety, like whether or not all points in a fibre of $\pi_B$ can be real, positive, or satisfy the simplicial analogue of the triangle inequality (i.e. are volumes which are \CHANGE{actually}{} realizable by a simplex). We take an experimental approach toward these questions for bases of $\mathcal M_3$.

\section{Notation and Background}
\label{sec:background}

\subsection{Distance geometry}
Given a collection $\mydef{\p}=(p_1,\ldots,p_{n+1})$
 of points in $\mathbb{R}^n$, not all on a hyperplane, their convex hull is an $n$-simplex which we denote by $\mydef{\Delta_n(\p)}\subseteq \mathbb{R}^n$. When discussing the $n$-simplex as an abstract polytope, or if $\p$ is understood to be arbitrary, we simply write $\mydef{\Delta_n}$. Faces of $\Delta_n$ are identified with subsets $S \subseteq \mydef{[n+1]}=\{1,2,\ldots,n+1\}$ and correspond to the simplices $\Delta_{|S|-1}(\p_S)$ obtained as convex hulls of the vertices $\mydef{\p_S}=(p_s)_{s \in S}$. 

We write the Euclidean volume of the face of $\Delta_n(\p) \subseteq \mathbb{R}^n$ indexed by $S$ as $\mydef{v_S}$. When explicitly listing $S$, we will often consolidate notation like $v_{123}\leftrightarrow v_{\{1,2,3\}}$ provided it does not introduce any ambiguity. In that case, the volumes are interpreted as indeterminants. Moreover, we write $\mydef{\x_S}=v_S^2$ for the corresponding squared volumes.   We will also use the notation $v_S$ and $x_S$ when referring to the volumes of faces of an arbitrary $n$-simplex.  Equipped with this notation, we state the eponymous result inspiring this paper.
\begin{theorem}[Heron \cite{Heron}]
\label{thm:Heron}
The squared area of any triangle $\Delta_2(p_1,p_2,p_3)$ may be expressed as
\[
x_{123}=\frac{1}{16}\cdot (2x_{12}x_{13}+2x_{12}x_{23}+2x_{13}x_{23}-x_{12}^2-x_{13}^2-x_{23}^2).
\]
\end{theorem}
Theorem \ref{thm:Heron} is referred to as \textit{Heron's formula}, after Heron of Alexandria (circa AD $62$) who wrote it in his book \textit{Metrica} \cite{Heron}. The following more general result is well-known.
\begin{theorem}[c.f. \cite{Blumenthal}]
\label{thm:CayleyMenger}
For any $n$-simplex $\Delta_n(\p)$, we have 
\[
x_{12\cdots n+1} = \frac{(-1)^{n+1}}{(n!)^2\cdot 2^n} \cdot  \begin{vmatrix}
0 & 1 & 1 & 1 & \cdots & 1 \\
1 & 0 & x_{1,2} & x_{1,3} & \cdots & x_{1,n+1} \\
1 & x_{1,2} & 0 & x_{2,3} & \cdots & x_{2,n+1} \\
\vdots & \vdots & \vdots & \vdots & \ddots & \vdots \\
1 & x_{1,n+1} & x_{2,n+1} & x_{3,n+1}  & \cdots & 0
\end{vmatrix}.
\]
\end{theorem}
The part of the matrix in Theorem \ref{thm:CayleyMenger} other than the first row and column is a \mydef{Euclidean distance matrix}. The whole matrix is called a \mydef{Cayley-Menger matrix} and its determinant is called a \mydef{Cayley-Menger determinant}.

Not every collection of ${{n+1}\choose{2}}$ positive real numbers are the edge lengths of a simplex. In $\mathbb{R}^2$ the condition on edge lengths to be those of a triangle \CHANGE{}{is known} as the \textit{triangle inequality}. In higher dimensions, this is known as Schoenberg's problem, whose solution comes in the form of classifying when a matrix is a Euclidean distance matrix.

\begin{theorem}[Shoenberg \cite{Shoenberg} (c.f. \cite{DistanceGeometry})]
\label{thm:Shoenberg}
The values $\{x_{ij}\}_{1 \leq i<j\leq n+1}$ are the squared edge lengths of a simplex if and only if the matrix $\frac{1}{2}(x_{1i}+x_{1j}-x_{ij} \mid 2 \leq i,j \leq n+1)$ is positive-definite of full rank.
\end{theorem}

We remark that Theorem \ref{thm:CayleyMenger} gives the squared volumes of all faces of $\Delta_n(\p)$ in terms of the principal minors of the Cayley-Menger matrix: the squared volume $x_{S}$ is, up to a constant, the principal minor of the matrix in Theorem \ref{thm:CayleyMenger} indexed by $\{0\} \cup S$. From this observation, we obtain a parametrization via these scaled principal minors:
\begin{align*}
\mydef{\varphi_n}:\mathbb{C}^{{n+1}\choose{2}} &\to \mathbb{C}^{2^{n+1}-(n+1)-1} \\ 
(x_{12},\ldots,x_{n,n+1}) &\mapsto (x_S)_{S \subseteq [n+1],|S|>1}.
\end{align*}
To simplify notation, let $\mydef{e(n)}={{n+1}\choose{2}}$ and $\mydef{N(n)} = 2^{n+1}-(n+1)-1$. Let $\mydef{\mathcal E_n} \subseteq \mathbb{R}_{\geq 0}^{e(n)}$ be the subset of values which satisfy Theorem \ref{thm:Shoenberg}. The \mydef{volume-realization space} of an $n$-simplex is the set $\mydef{X_n(\Delta)}=\varphi_n(\mathcal E_n)$. Similarly, we define $\mydef{X_n(\mathbb{R})}, \mydef{X_n(\mathbb{R}_{>0})},$ and $\mydef{X_n(\mathbb{C})}$ to be the images of $\mathbb{R}^{e(n)}$, $\mathbb{R}_{>0}^{e(n)}$, and $\mathbb{C}^{e(n)}$ under $\varphi_n$ respectively. We define $\mydef{X_n}$ to be the algebraic closure of $X_n(\mathbb{C})$, called the \mydef{$n$-th Heron variety}. We have the obvious set containments
\[
X_n(\Delta) \,\,\subsetneq\,\,  X_n(\mathbb{R}_{>0}) \,\,\subsetneq\,\, X_n(\mathbb{R}) \,\,\subsetneq\,\, X_n(\mathbb{C}) \,\,\subsetneq\,\, X_n \,\,\subsetneq\,\,  \mathbb{C}^{N(n)}.
\]
\begin{lemma}
\label{lem:HeronVariety}
The set $X_n$ is the Zariski closure of $X_n(\Delta)$.
\end{lemma}
\begin{proof}
$X_n$ is a graph over $\mathbb{C}^{e(n)}$ and thus has \CHANGE{}{(complex)} dimension $e(n)$ and is irreducible. The all ones vector $\textbf{1}$ belongs to the interior of $\mathcal E_n$  and so there is a\CHANGE{n}{ (real)} $e(n)$-dimensional neighborhood of $\textbf{1}$ in $\mathcal E_n$ whose image under $\varphi_n$ is thus also $e(n)$-dimensional. We have shown that $X_n(\Delta)$ contains a\CHANGE{n}{ (real)} $e(n)$-dimensional set and is contained in the \CHANGE{}{(complex)} $e(n)$-dimensional variety $X_n$ and so $\overline{X_n(\Delta)} = X_n$. 
\end{proof}

\begin{example}
The second Heron variety is the affine variety $X_2 \subseteq \mathbb{C}_{\textbf{x}}^4$ defined by Heron's formula
\[f_{123}(\textbf{x}) = 16x_{123}-2x_{12}x_{13}-2x_{12}x_{23}-2x_{13}x_{23}+x_{12}^2+x_{13}^2+x_{23}^2.
\]
It has degree $2$ and dimension $3$. 

The third Heron variety $X_3 \subseteq \mathbb{C}_{\textbf{x}}^{11}$ is cut out by five polynomials
{\footnotesize{
\begin{align*}
f_{123}(\textbf{x}) =& 16x_{123}-2x_{12}x_{13}-2x_{12}x_{23}-2x_{13}x_{23}+x_{12}^2+x_{13}^2+x_{23}^2 \\
f_{124}(\textbf{x}) =& 16x_{124}-2x_{12}x_{14}-2x_{12}x_{24}-2x_{14}x_{24}+x_{12}^2+x_{14}^2+x_{24}^2 \\
f_{134}(\textbf{x}) =& 16x_{134}-2x_{13}x_{14}-2x_{13}x_{34}-2x_{14}x_{34}+x_{13}^2+x_{14}^2+x_{34}^2 \\
f_{234}(\textbf{x}) =& 16x_{234}-2x_{23}x_{24}-2x_{23}x_{34}-2x_{24}x_{34}+x_{23}^2+x_{24}^2+x_{34}^2 \\
f_{1234}(\textbf{x}) = &288x_{1234}+2x_{12}x_{13}x_{23}-2x_{12}x_{14}x_{23}-2x_{13}x_{14}x_{23}+2x_{14}^2x_{23}+2x_{14}x_{23}^2\\ 
&-2x_{12}x_{13}x_{24}+2x_{13}^2x_{24}+2x_{12}x_{14}x_{24}-2x_{13}x_{14}x_{24}-2x_{13}x_{23}x_{24}\\
&-2x_{14}x_{23}x_{24}+2x_{13}x_{24}^2+2x_{12}^2x_{34}-2x_{12}x_{13}x_{34}-2x_{12}x_{14}x_{34}\\
&+2x_{13}x_{14}x_{34}-2x_{12}x_{23}x_{34}-2x_{14}x_{23}x_{34}-2x_{12}x_{24}x_{34}-2x_{13}x_{24}x_{34}\\
&+2x_{23}x_{24}x_{34}+2x_{12}x_{34}^2
\end{align*}
}}
It has dimension $6$ and degree $16$. 
\end{example}

\subsection{Algebraic matroids}
We streamline the necessary background on matroids and algebraic matroids in this section. A standard reference for the vast world of matroids is \cite{Oxley}. 
\begin{definition}[Matroid]
\label{def:matroid}
A \mydef{matroid} is a pair $\mathcal M=(E,\mathcal B)$ where $\mydef{E}$ is a finite set, called the \mydef{ground set} of $\mathcal M$, and $\mydef{\mathcal B}$ is a collection of subsets of $E$ called the \mydef{bases} of $\mathcal M$ which satisfy
\begin{enumerate}
\item $\mathcal B \neq \emptyset$
\item If $B_1,B_2 \in \mathcal B$ are distinct, then for any $a \in B_1\backslash B_2$, there exists $b \in B_2 \backslash B_1$ such that 
\[(B_1\backslash \{a\})\cup \{b\} \in \mathcal B.
\]
\end{enumerate}
\end{definition}
Given a matroid $\mathcal M=(E,\mathcal B)$, we define several other subsets:  subsets of $\mathcal B$ are called \mydef{independent sets}, those subsets of $E$ which are not independent are said to be  \mydef{dependent}, and minimal dependent sets are called \mydef{circuits}. Given a subset $S \subseteq E$, the maximal cardinality of an independent subset of $S$ is called the \mydef{rank} of $S$, denoted $\mydef{\textrm{rk}(S)}$.

\begin{definition}[Representable matroid]
Fix a field $\mathbb{F}$ and a matrix $A \in \mathbb{F}^{k \times n}$. The matroid $\mydef{\mathcal M(A)}$ associated to $A$ is the matroid on the ground set of columns of $A$ whose bases are those subsets of columns which form a basis for the column space of $A$. Any matroid which can be written as $\mathcal M(A)$ for some $A$, up to a relabeling of the ground set, is called \mydef{representable}.
\end{definition} We point-out that $\mathcal M(A)$ depends only on the row space of $A$, so  $\mydef{\mathcal M(\textrm{row}(A))}=\mathcal M(A)$ is well-defined.
 Representable matroids constitute only a small portion of those pairs satisfying Definition~\ref{def:matroid} \cite{BK,MNWW,Nelson}. Matroids can be constructed from a multitude of other mathematical objects, and are afforded names accordingly (e.g. graphic matroids, algebraic matroids, etc). 
Our focus is on \textit{algebraic matroids}. We define algebraic matroids over $\mathbb{C}$ geometrically. 
\begin{definition}
\label{def:algebraicmatroid}
Let $X \subseteq \mathbb{C}^{N}$ be an irreducible affine variety. The \mydef{algebraic matroid}  of $X$ is the pair $\mydef{\mathcal M(X)}=(E,\mathcal B)$, where $E=[N]$ are the indices of coordinates of its ambient space and $\mathcal B$ is the collection of subsets $B \subseteq E$ for which the coordinate projection
\begin{align*}
\mydef{\pi_{B}}:X &\to \mathbb{C}_{B}^{|B|}  \\
\textbf{x} &\mapsto  (x_b)_{b \in B}
\end{align*}
is finite-to-one and dominant.
\end{definition}
One may draw several immediate conclusions about an algebraic matroid $\mathcal M(X)=(E,\mathcal B)$ from Definition \ref{def:algebraicmatroid} along with the definition of a matroid. For example, all bases of $\mathcal M(X)$ have cardinality equal to the dimension of $X$, which is thus the rank of $\mathcal M(X)$. The rank of a subset $S$ of $E$ is the dimension of the image of $\pi_S(X)$, and $S$ is independent if $\pi_S$ is dominant. 

The main objects of interest in this paper are the algebraic matroids $\mydef{\mathcal M_n} = \mathcal M(X_n)$ coming from the Heron varieties.  The following fact provides the basis of our computation of $\mathcal M_n$.
\begin{lemma}{\cite[Proposition~6.7.10]{Oxley}}
\label{lem:tangentspace}
Let $X \subseteq \mathbb{C}^N$ be an irreducible affine variety. The matroid $\mathcal M(X)$ is the matroid of the tangent space of $X$ at a generic point. 
\end{lemma}
Since $X_n$ is parametrized by $\mathbb{C}^{e(n)}$ via the map $\varphi_n$, it is easy to write down the tangent space of $X_n$ at some point $\varphi(\textbf{x})$: it is the column space of the evaluated Jacobian matrix 
\[
\mathrm{Jac}(\varphi_n)|_{\textbf{x}} \in \mathbb{C}^{\CHANGE{e(n) \times N(n)}{N(n) \times e(n)}}.
\]
Hence, to compute the algebraic matroid $\mathcal M_n$ of the $n$-th Heron variety, one needs only to establish the algebraic matroid associated to (the transpose of) this Jacobian matrix. Since $X_n$ has dimension $e(n)$, so does its tangent space at a generic point, and therefore, the Jacobian of $\varphi_n$ is of full-rank. 
\begin{example}[The algebraic matroid of $X_2$]
\label{ex:X2Matroid}
Consider the (transpose) of the Jacobian of the map $\varphi_2:\mathbb{C}^3 \to \mathbb{C}^4$ parametrizing $X_2$:
\[
\textrm{Jac}(\varphi_2)^T = \begin{pmatrix}
1 & 0 & 0& (1/8)\cdot(-x_{12}+x_{13}+x_{23})  \\
0 & 1 & 0 &(1/8)\cdot(\,\,\,\,\,x_{12}-x_{13}+x_{23}) \\
 0 & 0 & 1&(1/8)\cdot(\,\,\,\,\,x_{12}+x_{13}-x_{23}) 
\end{pmatrix}.
\]
For generic values of $x_{12},x_{13},x_{23}$, every maximal minor of this matrix is nonzero, so by Lemma \ref{lem:tangentspace}, each $3$-subset of $E=\{12,13,23,123\}$ is a basis of the algebraic matroid $\mathcal M_2=\mathcal M(X_2)$:
\[\mathcal B=\{\{12,13,23\},\{12,13,123\},\{12,23,123\},\{13,23,123\}\}.\] 
The matroid $\mathcal M_2$ is also known as the \textit{uniform matroid} of rank $3$ on $4$ elements, $U(3,4)$. We point out that the (non-generic) values of $x_{12},x_{13},x_{23}$ for which the tangent space $T_{\varphi(\textbf{x})}(X_n)$ is not the uniform matroid are those squared edge lengths which correspond to right triangles. To see this, note that $-x_{12}+x_{13}+x_{23}\CHANGE{}{=0}$ is just \CHANGE{}{the} Pythagorean\CHANGE{'s}{} theorem with $x_{12}$ as the squared hypotenuse.
\end{example}

\begin{example}[The nontrivial basis of $\mathcal M_2$]
\label{ex:M2Basis}
Consider the basis $B=\{12,13,123\}$. The projection map $\pi_B:X_2 \to \mathbb{C}^3$ is  finite-to-one and dominant. Given a generic point $\textbf{a} = (a_{12},a_{13},a_{123}) \in \mathbb{C}^3$, the fibre  $\pi_B^{-1}(\textbf{a})$ consists of the two points of the form $(a_{12},a_{13},x_{23},a_{123})$ where 
\[
x_{23} = a_{12}+a_{13} \pm 2\sqrt{a_{12}a_{13}-4a_{123}}.
\]
Consequently $d_{B}=2$. Observe that the expression under the square root is zero whenever the corresponding triangle is a right triangle.
\end{example}

\begin{example} 
\label{ex:icecream}
Turning toward the matroid $\mathcal M_3$, consider the subset $B=\{12,13,14,123,124,134\}$ of the ground set $E=\{12,13,\ldots,1234\}$ of $\mathcal M_3$. This subset is illustrated in \CHANGE{Figure \ref{fig:icecream}}{Figure \ref{fig:icecream}}.
\begin{figure}[!htpb]
\includegraphics[scale=0.7]{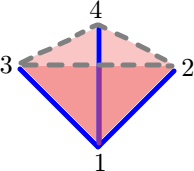}
\caption{A set $\{12,13,14,123,124,134\}$ of six faces of the tetrahedron $\Delta_3$.}
\label{fig:icecream}
\end{figure}

Since $\dim(X_3) = e(3) = 6$, the set $B$ has the necessary cardinality to be a basis of $\mathcal M_3$. Determining whether $B$ is indeed a basis involves computing the $6 \times 6$ determinant of  the submatrix of $\textrm{Jac}(\varphi_3)^T$ indexed by $B$:
\[{\small{
\begin{blockarray}{ccccccc}
&x_{12} & x_{13} & x_{14} & x_{123} & x_{124} & x_{134} \\
\begin{block}{c(cccccc)}
x_{12} & 1 & 0 & 0 & (1/8)\cdot(-x_{12}+x_{13}+x_{23}) & (1/8)\cdot(-x_{12}+x_{14}+x_{24}) & 0 \\
x_{13} & 0 & 1 & 0 & (1/8)\cdot(\,\,\,\,\,x_{12}-x_{13}+x_{23}) & 0 & (1/8)\cdot(-x_{13}+x_{14}+x_{34}) \\
x_{23} & 0 & 0 & 0 & (1/8)\cdot(\,\,\,\,\,x_{12}+x_{13}-x_{23}) & 0 & 0 \\
x_{14} & 0 & 0 & 1 & 0 & (1/8)\cdot(\,\,\,\,\,x_{12}-x_{14}+x_{24}) & (1/8)\cdot(\,\,\,\,\,x_{13}-x_{14}+x_{34}) \\
x_{24} & 0 & 0 & 0 & 0 & (1/8)\cdot(\,\,\,\,\,x_{12}+x_{14}-x_{24}) & 0 \\
x_{34} & 0 & 0 & 0 & 0 & 0 & (1/8)\cdot(\,\,\,\,\,x_{13}+x_{14}-x_{34})\\
\end{block}
\end{blockarray}}}
\]
This determinant is $(-1/8^3)(x_{13}+x_{14}-x_{34})(x_{12}+x_{14}-x_{24})(x_{12}+x_{13}-x_{23})$. Since this is generically nonzero,  $B$ is a basis. 
Notice that this minor of the Jacobian factors into three copies of the minor observed in Example \CHANGE{\ref{ex:M2Basis}}{\ref{ex:X2Matroid}}. This phenomenon is a consequence of the decomposability of $\pi_B$ as explored in Section \ref{secsec:decompose}.
\end{example}

\subsection{Branched covers}
Fix an irreducible affine variety $X \subset \mathbb{C}^n$ and a basis $B$ of the algebraic matroid of $X$. 
The projection maps $\pi_B:X \to \mathbb{C}_{B}^{|B|}$  are the main maps of interest for this work. They are finite-to-one dominant maps of irreducible affine varieties, also known as \mydef{branched covers}. Every branched cover represents an  \textit{enumerative problem}, which asks \textit{How many points are in  a generic fibre of $\pi_B$?} This number $\mydef{d_B}$ is the \mydef{degree} of $\pi_B$.

In the context of this subsection, we call the codomain $\mathbb{C}_{B}^{|B|}$ of $\pi_B$ the \mydef{parameter space} and its elements \mydef{parameters}.  The locus $\mydef{\textrm{Disc}} \subseteq \mathbb{C}_{B}^{|B|}$  of parameter values $\textbf{a} \in \mathbb{C}_{B}^{|B|}$  for which the fibre $\pi_B^{-1}(\textbf{a})$ does \textit{not} consist of $d_B$  isolated solutions is called the \mydef{discriminant} of $\pi_B$. The \mydef{branch locus} of $\pi_B$ is the subset of  $X$ for which the Jacobian of $\pi_B$ drops rank. The image of the branch locus is the discriminant.

\begin{example}[A decomposable branched cover]
Consider the variety parametrized by the map $(x,y) \mapsto (x,y,-x^4-yx^2)$. This is a hypersurface $X$ in $\mathbb{C}_{x,y,z}^3$ cut out by the single equation $x^4+yx^2+z=0$. The variables $B=\{y,z\}$ form a basis of $\mathcal M(X)$ and the map $\pi_B$ is clearly $4$-to-$1$. The branch locus is the intersection of $X$ with the locus for which the matrix
\[\begin{pmatrix} 0 & -4x^3-2xy \\ 1 & -x^2 \end{pmatrix}
\]
drops rank. Namely, the branch locus is $X \cap (\mathcal V(x) \cup \mathcal V(-4x^2-\CHANGE{y}{2y}))$. The projection of this branch locus onto $\mathbb{C}_{y,z}^2$ produces the discriminant $\textrm{Disc} = \mathcal V(z) \cup \mathcal V(y^2-4z)$.

Over any parameter value not in the discriminant, there are four solutions which can be written in terms of radicals by first solving for $w=x^2$ and then taking roots:
\[
x = \pm \sqrt{w} = \pm \sqrt{\frac{-y \pm \sqrt{y^2-4z}}{2}}.
\]
The substitution $w=x^2$ may be thought of as a decomposition of $\pi_B$ as 
\[
\begin{tikzpicture}
\node at (0,0) (a) {$\mathcal V(x^4+yx^2+z)$};
\node[right =1cm of a]  (b){$\mathcal V(w^2+yw+z)$};
\node[right =1cm of b]  (c){$\mathbb{C}_{y,z}^2$.};
\draw[->,>=stealth] (a) --node[above]{$\psi$} (b);
\draw[->,>=stealth] (a) --node[below]{$2:1$} (b);
\draw[->,>=stealth] (b) --node[above]{$\phi$} (c);
\draw[->,>=stealth] (b) --node[below]{$2:1$} (c);
\draw[->,>=stealth] (a) edge[bend right=40]node[above]{$\pi_B$}node[below]{$4:1$} (c);
\end{tikzpicture}
\]Branched covers which decompose into compositions of branched covers, each with degree $>1$, are called \textit{decomposable}. For more details on decomposable branched covers, see \cite{SparseDecomposable}.
As a consequence of this decomposition, the Jacobian of $\pi_B$ factors into the product of the Jacobians of $\psi$ and $\phi$. The two components $\mathcal V(z), \mathcal V(y^2-4z)$ of the discriminant can thus be attributed to the loci where $\textrm{Jac}(\psi)$ and $\textrm{Jac}(\phi)$ drop rank, respectively. 
\end{example}

The complement of $\textrm{Disc}$ is the set $\mydef{\mathcal U}$ of regular values of $\pi_B$. The set $\mathcal U$ is path-connected over $\mathbb{C}^{|B|} \cong \mathbb{R}^{2|B|}$ since the discriminant is at most a \textit{complex} hypersurface. Locally, the map $\pi_B$ restricted to the preimage of $\mathcal U$ is a $d_B$-to-one \textit{covering space}.  

\subsection{Permutation groups}
We recall some background on permutation groups (see \cite{PermBackground} for details). A \mydef{permutation group} is any subgroup $G$ of a symmetric group $\mydef{\mathfrak S_d}$ on $d$ objects. We assume these objects are labeled by $[d]$. The action of $\mathfrak S_d$ on $[d]$ descends to an action of $G$ on $[d]$ via the inclusion $G \hookrightarrow \mathfrak S_d$. The group $G$ is said to be \mydef{transitive} if there is only one orbit of this group action. 

\CHANGE{}{Assume $G\subseteq \mathfrak S_d$ is transitive.}
Any set partition $P \CHANGE{=}{: } P_1 \sqcup P_2 \cdots \sqcup P_k$ of $[d]$ is called a \mydef{block system} of $G$ if the partition $P$ is \mydef{$G$-invariant}: for all $g \in G$, we have $x,y \in P_i \iff gx,gy \in P_j$ for some $j$. The parts $P_1,\ldots,P_k$ are called \mydef{blocks}, and the action of $G$ on $[d]$ extends to an action of $G$ on $\{P_1,\ldots,P_k\}$. If the former action is transitive, then clearly so is the latter and it follows that all blocks in a block system have the same size. If $G$ has a nontrivial block system (comprised of more than one block of size larger than one) then $G$ is said to be \mydef{imprimitive} and otherwise it is \mydef{primitive}. Given a partition $P$ of $[d]$, we write $\mydef{\mathfrak S_P} \subseteq \mathfrak S_{d}$ for the subgroup of permutations which preserve $P$. 

\subsection{Monodromy}
Consider an affine variety $X$ and a basis $B$ of its algebraic matroid $\mathcal M(X)$. The branched cover $\pi_B: X \to \mathbb{C}_{B}^{|B|}$, restricted to the preimage of its regular values $\mathcal U \subseteq \mathbb{C}_{B}^{|B|}$, is a $d_B$-to-one covering space. As a consequence, any path $\gamma:[0,1]_t \to \mathcal U$ lifts to $d_B$ well-defined paths $\textbf{x}^{(1)}(t),\ldots,\textbf{x}^{(d_B)}(t)$ via analytic continuation. If $\gamma$ is a loop based at a point $\textbf{u} \in \mathcal U$, then the map $\textbf{x}^{(i)}(0)\mapsto \textbf{x}^{(i)}(1) = \textbf{x}^{(j)}(0)$ is a permutation of the fibre $\pi_B^{-1}(\textbf{u})$. After labeling the points of this fibre like $i \leftrightarrow \textbf{x}^{(i)}(0)$, we write the permutation induced by $\gamma$ as  $\sigma_\gamma \in\mathfrak S_{d_B}$. The collection of all permutations of the form $\sigma_\gamma$ for some path $\gamma:[0,1]_t \to \mathcal U$  based at $\textbf{u}$ is called the \mydef{monodromy group of $\pi_B$ based at $\textbf{u}$}. This group is well-defined as a subgroup of $\mathfrak S_{d_B}$ up to relabeling of the points in the fibre. In fact, up to relabeling, this group does not depend on the base point at all. We call this permutation group \mydef{$G_{\pi_B}$}the \mydef{monodromy group} of $\pi_B$. 

The monodromy group of $G_{\pi}$ is transitive since the preimage of $\mathcal U$ is path-connected since $X$ is irreducible: to connect two points in a fibre via a monodromy permutation $\sigma_\gamma$, simply choose the loop $\gamma$ which is the image of any \CHANGE{real}{}path $\eta:[0,1] \to X$ connecting them. 

\subsection{Decomposability}
\label{secsec:decompose}
The dictionary between geometric properties, like the irreducibility of $X$, and group-theoretic properties, like transitivity of $G_{\pi}$, runs deep. The next result describes the connection between \textit{decomposability} and \textit{imprimitivity}. Its proof is straightforward and omitted. A branched cover $\pi$ is \mydef{decomposable} when $\pi=\phi \circ \psi$ for some branched covers $\phi$ and $\psi$, of degree $>1$.
\begin{lemma}
\label{lem:blocksfromdecompositions}
Suppose $\pi=\phi \circ \psi$ is a branched cover and $G_{\pi} \subset \mathfrak S_{\pi^{-1}(\textbf{u})} \cong \mathfrak S_{d}$ for some regular value $\textbf{u}$ of $\pi$. Then  $P_{\phi,\psi}=\bigsqcup_{y \in \phi^{-1}(\textbf{u})} \psi^{-1}(y)$ is a block system of \CHANGE{$\pi$}{$G_\pi$}. In particular, if $\pi$ is decomposable, then $G_{\pi}$ is imprimitive. 
\end{lemma}
The converse of the last part of Lemma \ref{lem:blocksfromdecompositions}, that imprimitivity \textit{implies} decomposability, is true as well \cite{SparseDecomposable, Pirola}. Combining Lemma~\ref{lem:blocksfromdecompositions} with the observation that every coordinate projection $\pi_B:X \to \mathbb{C}_{B}^{|B|}$ corresponding to a basis $B=\{b_1,\ldots,b_r\}$ is obtained as a composition of projections 
\[
\begin{tikzpicture}
\node at (0,0) (a) {$X$};
\node[right =1cm of a]  (b){$X^{(1)}$};
\node[right =1cm of b]  (c){$X^{(1,2)}$};
\node[right =1cm of c]  (d){$\cdots$};
\node[right =1cm of d]  (e){$X^{(1,2,\ldots,N-r)}=\mathbb{C}_{B}^{|B|}$};
\draw[->,>=stealth] (a) --node[above]{} (b);
\draw[->,>=stealth] (a) --node[below]{} (b);
\draw[->,>=stealth] (b) --node[above]{} (c);
\draw[->,>=stealth] (b) --node[below]{} (c);
\draw[->,>=stealth] (c) --node[below]{} (d);
\draw[->,>=stealth] (d) --node[below]{} (e);
\draw[->,>=stealth] (a) edge[bend right=20]node[above]{$\pi_B$}node[below]{$d_B:1$} (e);
\end{tikzpicture}
\]
yields the following corollary.
\begin{corollary}
\label{cor:coordinatesymmetries}
Let $S$ be a generic fibre of $\pi_B$ and for each $i \in [n]$ let $P_i$ be the partition of $S$ induced by grouping elements of $S$ together if they agree in their $i$-th coordinate. Then each $P_i$ is a block system for $G_{\pi_B}$ and so
\[
G_{\pi_B} \subseteq \bigcap_{i=1}^n \mathfrak S_{P_i}=: \mydef{\widehat{G}_{\pi_B}}.
\]
\end{corollary}
We call the group $\widehat{G}_{\pi_B}$ in Corollary \ref{cor:coordinatesymmetries} the \mydef{coordinate symmetry group} of $\pi_B$. 

\begin{remark}
Corollary \ref{cor:coordinatesymmetries} can be easily extended by replacing the phrase \textit{agree in their $i$-th coordinate} with \textit{have the same value under a rational function $f:X \to \mathbb{C}$}. For our purposes, we only consider the rational functions given by the coordinate functions on $X$. 
\end{remark}

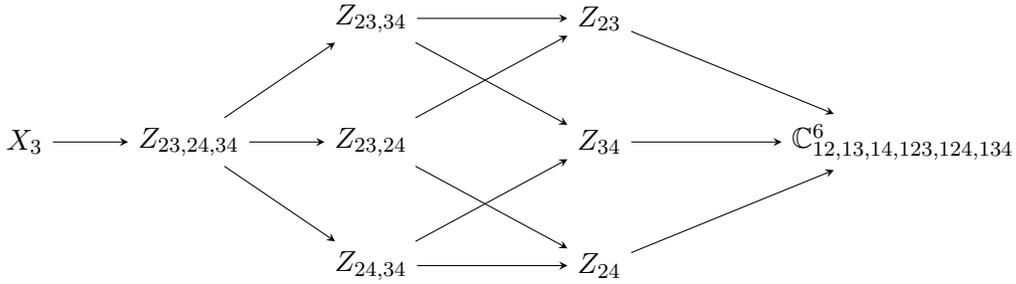
\begin{figure}[!htpb]
\[
\begin{tikzpicture}
\node at (-1,0) (q) {$X_3$};
\node[right=1cm of q] (a){$Z_{23,24,34}$};
\node[right =1cm of a]  (b){$Z_{23,24}$};
\node[above =1cm of b]  (c){$Z_{23,34}$};
\node[below =1cm of b]  (d){$Z_{24,34}$};
\node[right =2cm of c]  (e){$Z_{23}$};
\node[right =2cm of d]  (f){$Z_{24}$};
\node[right =2cm of b]  (g){$Z_{34}$};
\node[right =2cm of g]  (h){$\mathbb{C}_{12,13,14,123,124,134}^6$};
\draw[->,>=stealth] (a) -- (b);
\draw[->,>=stealth] (q) -- (a);
\draw[->,>=stealth] (a) -- (c);
\draw[->,>=stealth] (a) -- (d);
\draw[->,>=stealth] (b) -- (e);
\draw[->,>=stealth] (b) -- (f);
\draw[->,>=stealth] (c) -- (g);
\draw[->,>=stealth] (c) -- (e);
\draw[->,>=stealth] (d) -- (g);
\draw[->,>=stealth] (d) -- (f);
\draw[->,>=stealth] (e) -- (h);
\draw[->,>=stealth] (f) -- (h);
\draw[->,>=stealth] (g) -- (h);
\end{tikzpicture}
\]
\caption{The structure of several non-trivial decompositions of the branched cover $\pi_B:X_3 \to \mathbb{C}_B^6$ associated to the basis $B = \{12,13,14,123,124,134\}$ of $\mathcal M_3$. Each arrow, other than the left-most, represents a $2$-to-$1$ map.}
\label{fig:decomposedicecream}
\end{figure}

\begin{example}
\label{ex:icecreamrevisited}
We return to Example \ref{ex:icecream} regarding the basis $B=\{12,13,14,123,124,134\}$ of $\mathcal M_3$ and its corresponding projection $\pi_B$. Given generic values of the squared edge lengths $\{v_{12},v_{13},v_{14}\}$ and squared triangular areas $\{v_{123},v_{124},v_{134}\}$, the remaining three squared edge lengths may be found via Heron's formula as in Example \ref{ex:M2Basis}. 
Hence, the branched cover $\pi_B$ decomposes in several different ways according to Figure \ref{fig:decomposedicecream}, where $Z_{S}$ is the image of $X_3$ under the projection onto the variables $B \cup S$. All maps other than the leftmost one are $2$-to-$1$. Consequently, this figure shows that $G_{\pi_B}$ is imprimitive. 

Figure \ref{fig:icecreamBigDecomposition} refines Figure \ref{fig:decomposedicecream} by tracking which coordinates of the eight points in a generic fibre of $\pi_B$ share equal values: points labeled $s_{ijk}$ and $s_{i'j'k'}$ share a value in the $24,34,$ or $23$ coordinate if and only if $i=i', j=j',$ or $k=k'$ respectively. The automorphisms of the graph in Figure \ref{fig:icecreamBigDecomposition} is the group which preserves all of the block systems implied by the fibre structures  of the decompositions of $\pi_B$. Consequently, this diagram describes the coordinate \CHANGE{monodromy}{symmetry} group $\widehat{G}_{\pi_B}$.

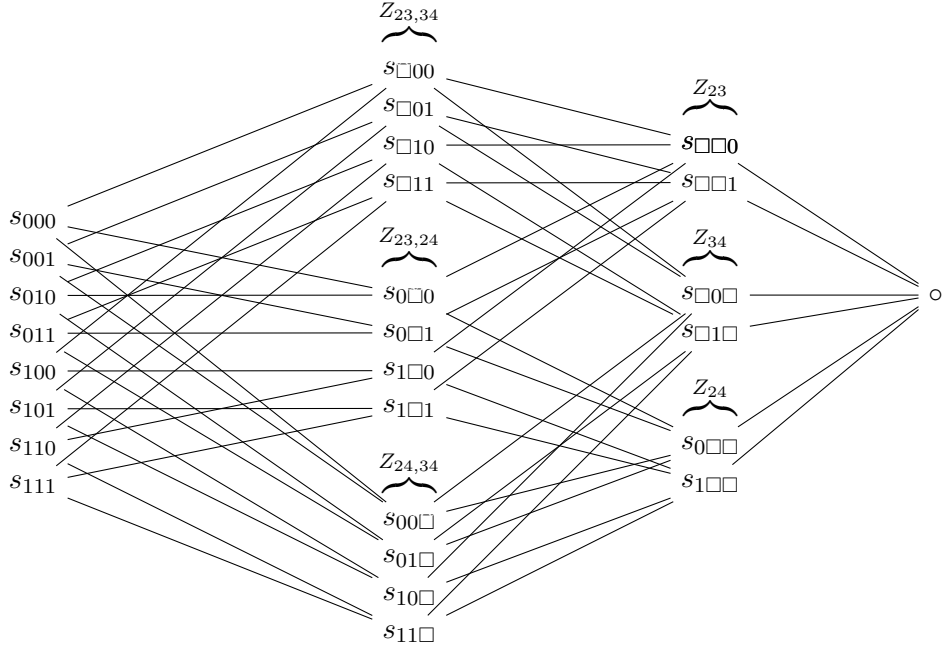
\begin{figure}[!htpb]
\begin{tikzpicture}
\node at (-1,1) (q000) {$s_{000}$};
\node[below=0.01cm of q000] (q001){$s_{001}$};
\node[below=0.01cm of q001] (q010){$s_{010}$};
\node[below=0.01cm of q010] (q011){$s_{011}$};
\node[below=0.01cm of q011] (q100){$s_{100}$};
\node[below=0.01cm of q100] (q101){$s_{101}$};
\node[below=0.01cm of q101] (q110){$s_{110}$};
\node[below=0.01cm of q110] (q111){$s_{111}$};
\node at (4,3) (qx00) {$s_{\square 0 0}$};
\node[below=0.01cm of qx00] (qx01){$s_{\square 0 1}$};
\node[below=0.01cm of qx01] (qx10){$s_{\square 1 0}$};
\node[below=0.01cm of qx10] (qx11){$s_{\square 1 1}$};
\node at (4,0) (q0x0) {$s_{0 \square 0}$};
\node[below=0.01cm of q0x0] (q0x1){$s_{0 \square 1}$};
\node[below=0.01cm of q0x1] (q1x0){$s_{1 \square 0}$};
\node[below=0.01cm of q1x0] (q1x1){$s_{1 \square 1}$};
\node at (4,-3) (q00x) {$s_{0 0\square}$};
\node[below=0.01cm of q00x] (q01x){$s_{0 1 \square}$};
\node[below=0.01cm of q01x] (q10x){$s_{1 0 \square}$};
\node[below=0.01cm of q10x] (q11x){$s_{1 1 \square}$};
\node at (8,2.5) (label23) {$\overbrace{{\color{white}{blah}}}^{Z_{23}}$};
\node at (8,-1.5) (label24) {$\overbrace{{\color{white}{blah}}}^{Z_{24}}$};
\node at (8,0.5) (label34) {$\overbrace{{\color{white}{blah}}}^{Z_{34}}$};
\node at (8,2) (qxx0) {$s_{\square \square 0}$};
\node at (4,3.5) (label2334) {$\overbrace{{\color{white}{blah}}}^{Z_{23,34}}$};
\node at (4,0.5) (label2324) {$\overbrace{{\color{white}{blah}}}^{Z_{23,24}}$};
\node at (4,-2.5) (label2434) {$\overbrace{{\color{white}{blah}}}^{Z_{24,34}}$};
\node at (8,2) (qxx0) {$s_{\square \square 0}$};
\node[below=0.01cm of qxx0] (qxx1){$s_{\square \square 1}$};
\node at (8,0) (qx0x) {$s_{\square  0 \square}$};
\node[below=0.01cm of qx0x] (qx1x){$s_{\square  1 \square}$};
\node at (8,-2) (q0xx) {$s_{0 \square \square} $};
\node at (11,0) (qxxx) {$\circ$};
\node[below=0.01cm of q0xx] (q1xx){$s_{1 \square \square}$};
\draw[-,>=stealth] (q000) -- (qx00);
\draw[-,>=stealth] (q000) -- (q0x0);
\draw[-,>=stealth] (q000) -- (q00x);
\draw[-,>=stealth] (q001) -- (qx01);
\draw[-,>=stealth] (q001) -- (q0x1);
\draw[-,>=stealth] (q001) -- (q00x);
\draw[-,>=stealth] (q010) -- (qx10);
\draw[-,>=stealth] (q010) -- (q0x0);
\draw[-,>=stealth] (q010) -- (q01x);
\draw[-,>=stealth] (q100) -- (qx00);
\draw[-,>=stealth] (q100) -- (q1x0);
\draw[-,>=stealth] (q100) -- (q10x);
\draw[-,>=stealth] (q011) -- (qx11);
\draw[-,>=stealth] (q011) -- (q0x1);
\draw[-,>=stealth] (q011) -- (q01x);
\draw[-,>=stealth] (q110) -- (qx10);
\draw[-,>=stealth] (q110) -- (q1x0);
\draw[-,>=stealth] (q110) -- (q11x);
\draw[-,>=stealth] (q111) -- (qx11);
\draw[-,>=stealth] (q111) -- (q1x1);
\draw[-,>=stealth] (q111) -- (q11x);
\draw[-,>=stealth] (q101) -- (qx01);
\draw[-,>=stealth] (q101) -- (q1x1);
\draw[-,>=stealth] (q101) -- (q10x);
\draw[-,>=stealth] (q00x) -- (q0xx);
\draw[-,>=stealth] (q00x) -- (qx0x);
\draw[-,>=stealth] (q01x) -- (q0xx);
\draw[-,>=stealth] (q01x) -- (qx1x);
\draw[-,>=stealth] (q10x) -- (q1xx);
\draw[-,>=stealth] (q10x) -- (qx0x);
\draw[-,>=stealth] (q11x) -- (q1xx);
\draw[-,>=stealth] (q11x) -- (qx1x);
\draw[-,>=stealth] (q0x0) -- (q0xx);
\draw[-,>=stealth] (q0x0) -- (qxx0);
\draw[-,>=stealth] (q0x1) -- (q0xx);
\draw[-,>=stealth] (q0x1) -- (qxx1);
\draw[-,>=stealth] (q1x0) -- (q1xx);
\draw[-,>=stealth] (q1x0) -- (qxx0);
\draw[-,>=stealth] (q1x1) -- (q1xx);
\draw[-,>=stealth] (q1x1) -- (qxx1);
\draw[-,>=stealth] (qx00) -- (qx0x);
\draw[-,>=stealth] (qx00) -- (qxx0);
\draw[-,>=stealth] (qx01) -- (qx0x);
\draw[-,>=stealth] (qx01) -- (qxx1);
\draw[-,>=stealth] (qx10) -- (qx1x);
\draw[-,>=stealth] (qx10) -- (qxx0);
\draw[-,>=stealth] (qx11) -- (qx1x);
\draw[-,>=stealth] (qx11) -- (qxx1);
\draw[-,>=stealth] (qxx0) -- (qxxx);
\draw[-,>=stealth] (qxx1) -- (qxxx);
\draw[-,>=stealth] (qx0x) -- (qxxx);
\draw[-,>=stealth] (qx1x) -- (qxxx);
\draw[-,>=stealth] (q0xx) -- (qxxx);
\draw[-,>=stealth] (q1xx) -- (qxxx);
\end{tikzpicture}
\caption{A refinement of Figure \ref{fig:decomposedicecream} describing how the fibres of the various decompositions of $\pi_B$ interact, where $B$ is the basis of $\mathcal M_3$ in Example \ref{ex:icecream}.}
\label{fig:icecreamBigDecomposition}
\end{figure} 

The permutations of the nodes on the left of the graph in Figure \ref{fig:icecreamBigDecomposition} induced by automorphisms of the graph establishes that the coordinate symmetry group of $\pi_B$ is the transitive permutation group $\widehat{G}_{\pi_B} = \ZZZ \times \ZZZ \times \ZZZ \subset \mathfrak S_8$.  We show that this equals the Galois group in Theorem \ref{thm:GaloisX3} using numerical algebraic geometry.
\end{example}

\begin{figure}[!htpb]
${\footnotesize{
\begin{blockarray}{cccccccccc}
& s_1 & s_2 & s_3 & s_4 & s_5 & s_6 & s_7 & s_8 \\
\begin{block}{c(ccccccccc)}
x_{23} & -1.5+0.8i &   -1.5+0.8i &   -1.5+0.8i &   -1.5+0.8i &    1.5-3.1i &    1.5-3.1i &   1.5-3.1i &    1.5-3.1i  \\
x_{24} & -0.5+2.0i &   -0.5+2.0i &    0.9-2.4i &    0.9-2.4i &   -0.5+2.0i &   -0.5+2.0i &   0.9-2.4i &    0.9-2.4i  \\
x_{34} &  5.0+0.1i &   -5.1+0.6i &    5.0+0.1i &   -5.1+0.6i &    5.0+0.1i &   -5.1+0.6i &   5.0+0.1i &   -5.1+0.6i \\
\end{block}
\end{blockarray}}}$

\includegraphics[scale=0.21]{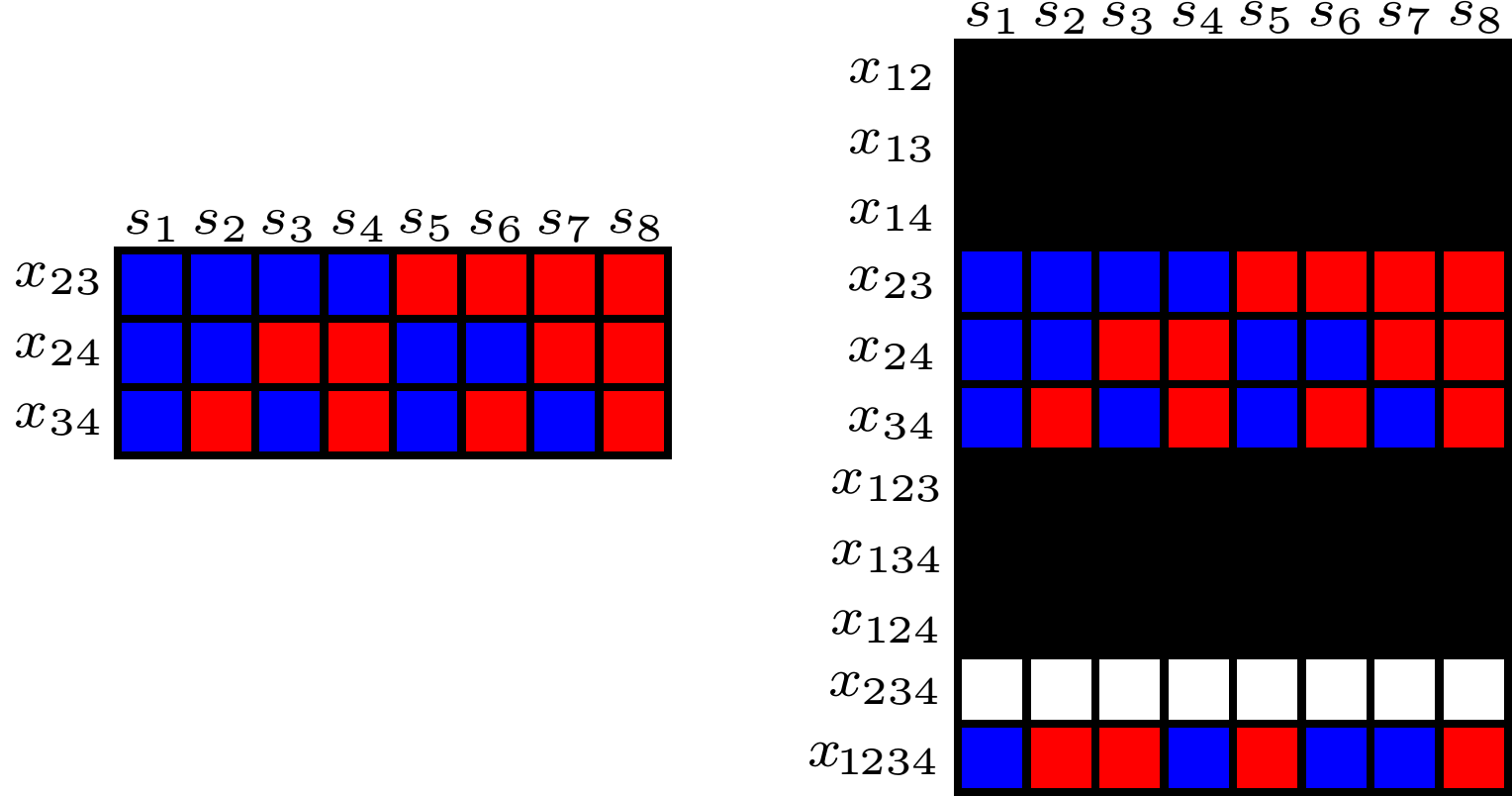}
\caption{(Top) A matrix showing some coordinates (rounded to fit on the page) of eight points in a generic fibre of $\pi_B$ for the basis  $B=\{12,13,23,123,124,134\}$  of $\mathcal M_3$. (Bottom Left) The symmetry diagram corresponding to those coordinates. (Bottom Right) The symmetry diagram corresponding to all coordinates; the white row indicates that the coordinate $x_{234}$ is distinct for each solution and the black rows indicate that those coordinates are equal amongst all solutions. }
\label{fig:coordinatesymmetry}
\end{figure}

Another way of illustrating the coordinate \CHANGE{monodromy}{symmetry} group of $\pi_B$ is by a diagrammatical array of colours indicating which coordinates of the points in a generic fibre agree. More specifically, consider the subset $\mydef{\hat{\textbf{x}}_B}$ of coordinates of $\mathbb{C}^{N(n)}$ for which points in a generic fibre of $\pi_B$ share a value.  Then write a $|\hat{\textbf{x}}_B| \times d_B$ grid whose columns correspond to points in a generic fibre and whose rows correspond to the coordinates in $\hat{\textbf{x}}_B$. Two boxes in the $x_i$-th row are coloured the same colour if the points corresponding to their columns agree in that coordinate. We call this the \mydef{coordinate symmetry diagram} of the branched cover $\pi_B$. For example, Figure \ref{fig:coordinatesymmetry} shows the coordinates which agree on a generic fibre of $\pi_B$ for $B = \{12,13,14,123,124,134\}$ along with  a coordinate symmetry diagram corresponding to the missing squared edge lengths $x_{23},x_{24},$ and $x_{34}$ (left) and one corresponding to all simplex volumes (right). For the coordinate symmetry diagrams associated to other bases, see Theorem \ref{thm:CoordinateSymmetryX3} and Figure \ref{fig:symmetrydiagrams}.

The monodromy group $G_{\pi}$ of a branched cover $\pi:X \to \mathbb{C}^r$ is a Galois group. Specifically, in \cite{Harris} Harris showed that $G_{\pi}$ is the Galois group of the field extension $\mathbb{C}(X)/\mathbb{C}(\mathbb{C}^r)$. As a consequence, a formula in radicals for \textit{all} coordinates of points in a generic fibre $\pi^{-1}(\textbf{u})$ exists if and only if the monodromy group $G_{\pi}$ is \textit{solvable}. See \cite{Zhang2019GaloisGO} for additional details and exposition regarding this correspondence.

\subsection{Numerical Algebraic Geometry}

Given an irreducible variety $X \subseteq \mathbb{C}^N$ of codimension $k$ defined by some polynomial system $F=(f_1,\ldots,f_k) \in \mathbb{C}[x_1,\ldots,x_N]$, a common way to represent the branched cover $\pi_B:X \to \mathbb{C}_{B}^{|B|}$ is by simply declaring the variables indexed by $B$ to be parameters. That is, by realizing the system $F$ as a $0$-dimensional \textit{parametrized polynomial system} in $\mathbb{C}[\textbf{x}_{B}][\textbf{x}_{[N]-B}]$.

As outlined in Section 3.2 of \cite{MetricAG}, a parametrized polynomial system like $F\subset \mathbb{C}[\textbf{x}_B][\textbf{x}_{[N]-B}]$ is the context in which one may apply a \textit{parameter homotopy} to follow numerical approximations of fibres of $\pi_B$ through the parameter space $\mathcal U$. Finding a single fibre can be done via a standard total-degree homotopy (see Section 3.3 of \cite{MetricAG}). These are standard tools in the world of \textit{numerical algebraic geometry}.  We use \texttt{HomotopyContinuation.jl} \cite{HCjl} through \texttt{julia} for all of our core numerical algebraic geometry computations. Specifically, this package is the numerical backbone of the new software system \texttt{Pandora.jl} \cite{Pandora} designed for studying enumerative problems. 

Equipped with the ability to follow the fibres of $\pi_B$ over the regular values $\mathcal U$, there is a straightforward algorithm for computing a monodromy permutation $\sigma_{\gamma}$ given a loop $\gamma:[0,1]_t \to \mathcal U$: simply follow the fibre $\pi^{-1}_{B}(\gamma(0))$ as $t$ goes from $0$ to $1$ along the $d_B$ solution paths $s_1(t),\ldots,s_{d_B}(t)$, and construct the permutation $s_i(0)\to s_i(1)$ (see Algorithm \ref{alg:NumericalMonodromyGroup}). 

\section{Computational Pipeline}
\label{sec:pipeline}
In this section we lay out a computational pipeline for addressing questions $\textbf{Q1}-\textbf{Q3}$, combining techniques in computational group theory, computational algebraic geometry, polyhedral geometry, and numerical algebraic geometry.

\subsection{Step 1: Determining candidate basis orbits} The labeling of the vertices of an $n$-simplex is arbitrary for our purposes: the \CHANGE{an}{}action of $\mathfrak S_{n+1}$ on the indices of variables in $\mathbb{C}_{\textbf{x}}^{N(n)}$ leaves the $n$-th Heron variety $X_n$ invariant. Since the rank of $\mathcal M_n$ is $e(n) = {{n+1}\choose{2}}$, our goal is to establish which, of the $\mydef{\alpha_n}={{N(n)}\choose {e(n)}}$-many subsets $B$ of size $e(n)$, are bases of $\mathcal M_n$. Using \texttt{GAP} through the \texttt{julia} package \texttt{OSCAR}, we can reduce this number from $\alpha_n$ to $\mydef{\beta_n}$: the number of \textit{orbits} of subsets of size $e(n)$ under the action of $\mathfrak S_{n+1}$ on \textit{sets of subsets} of $[n+1]$.
\[
\begin{tabular}{|c||c|c|c|}\hline
$n$ & 2 & 3 & 4  \\ \hline
$\alpha_n$ & 4 & 462 & 5311735 \\
$\beta_n$ & 2 & 35 & 48533\\
$\beta_n/\alpha_n$ & 0.5 & 0.06 & 0.009\\ \hline
\end{tabular}
\]
For example, there are four $3$-subsets of the positive dimensional faces of $\Delta_2$: 
\[\{12,13,23\},\{12,13,123\},\{12,23,123\},\{13,23,123\},
\] but the last three are in the same orbit under the action of $\mathfrak S_3$ on the vertices of $\Delta_2$. 
Na\"ively applying a function like \texttt{minimal\_image} in \texttt{GAP} \cite{GAP} to compute minimal images of each of the $\alpha_n$ sets of subsets of $[n+1]$ of size $e(n)$ and subsequently keeping the unique ones works for $n=3$ but is burdensome for $n=4$.
Hence, we divide the problem based on the $f$-vector of a set of faces of $\Delta_n$. The \mydef{$f$-vector} of a set $B$ of faces is the vector $\mydef{f(B)} = (\mydef{f_1(B)},\ldots,\mydef{f_n(B)})$ where $f_i(B)$ is the number of faces in $B$ of dimension $i$.

\begin{example}
Let $n=3$ and consider all $f$-vectors of $6$-subsets of positive dimensional faces of $\Delta_3$. These are the partitions $\lambda=(\lambda_1,\lambda_2,\lambda_3)$ of $6$ whose $i$-th part $\lambda_i$ is no larger than the $i$-th part of $f(\Delta_3) = (6,4,1)$. There are $10$ such possible vectors, which partition the set of all $35$ sextuples of faces of $\Delta_3$:
\[{\footnotesize{
\begin{tabular}{|c||c|c|c|c|c|c|c|c|c|c|} \hline
$f$-vector & (6, 0, 0)& (5, 1, 0)& (5, 0, 1)& (4, 2, 0)& (2, 4, 0)& (4, 1, 1)& (1, 4, 1)& (3, 3, 0)& (3, 2, 1)& (2, 3, 1) \\ \hline
\# orbits & 1 & 2 & 1 & 6 & 2 & 4 & 1 & 6 & 8 & 4\\  \hline 
\end{tabular}}}
\]
\end{example}

The upgraded algorithm for finding orbit representatives of each $e(n)$-subset of the $N(n)$-many positive dimensional faces of $\Delta_n$ is given in Algorithm \ref{alg:OrbitsWithFvector}. The idea is to first compute all possible $f$-vectors of such subsets, and then subsequently compute all orbits with that $f$-vector by first reducing the computation via taking the orbits of some particular dimensional part first. The benefit is illustrated in the following example.

\begin{example}
Suppose we want to compute all orbit representatives of $6$-subsets of positive dimensional faces of $\Delta_3$ with $f$-vector $(3,2,1)$. Instead of first listing the $\binom{6}{3}\cdot \binom{4}{2} \cdot \binom{1}{1} = 20 \cdot 6 \cdot 1 = 120$ subsets with that $f$-vector, and subsequently taking minimal images with respect to the action of $\mathfrak S_4$ on the face lattice of $\Delta_3$, we choose a triple of edges \textit{first}. We list the $20$ triples of edges of $\Delta_3$, and take $20$ minimal images to obtain $3$ orbits \[O_{1}^{(1)}=|\{12,13,14\}|,\quad \quad \quad O_2^{(1)}=|\{12,23,34\}|,\quad \quad \quad O_3^{(1)}=|\{12,23,13\}|
\] of triples of edges. 
We then consider the $3 \cdot 6 \cdot 1=18$ many subsets of the form $O_{i}^{(1)} \cup F_2 \cup F_3$ where $F_2$ is a pair of triangular faces and $F_3$ is the only $3$-dimensional face of $\Delta_3$. We compute $18$ additional minimal images and remove all duplicates. In total, instead of performing $120$ minimal image calls, we only performed $38$. The $8$ orbits of $6$-subsets of positive dimensional faces of $\Delta_3$ with $f$-vector $(3,2,1)$ are shown in Figure \ref{fig:EdgeTriples} under the orbit of triples of edges they were derived from.

\begin{figure}[!htpb]
\includegraphics[scale=0.3]{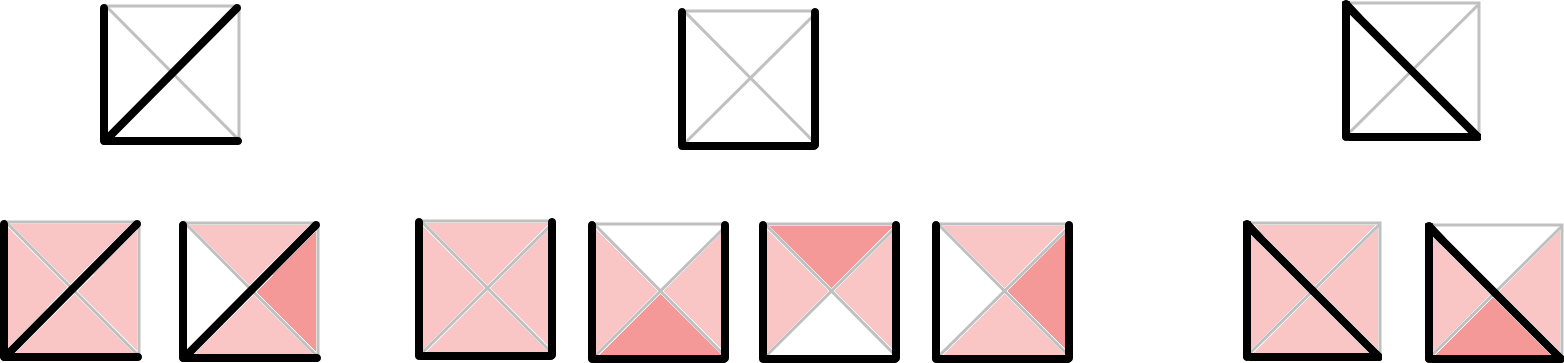}
\caption{Orbit representatives of the three orbits of triples of edges of $\Delta_3$ under the action of $\mathfrak S_{4}$ on the vertices of a tetrahedron.}
\label{fig:EdgeTriples}
\end{figure}
\end{example}

\begin{algorithm}[!htb]
\small
\caption{{\sf Orbits with $f$-vector}}\label{alg:OrbitsWithFvector}
\algnotext{EndIf}
\algnotext{EndFor}
\algnotext{EndWhile}
\begin{flushleft}
\textbf{Input: }{A partition $\lambda=(\lambda_1,\ldots,\lambda_n)$ of $e(n)$}  \\
\textbf{Output: }{An orbit representative for each orbit of $\mathfrak S_{n+1}$ acting on sets of faces of $\Delta_n$ with $f$-vector $\lambda$} \\
\end{flushleft}
\begin{algorithmic}[1]
\State Determine the $\lambda_i$ for which $\binom{f_i(\Delta_n)}{\lambda_i}$ is maximal 
\State Compute an $\mathfrak S_{n+1}$-orbit representative $O_S$ for each $\lambda_i$-subset $S$ of $i$-faces of $\Delta_n$
\State Define $\texttt{Orb}_i$ to be the set of all representatives \Comment{i.e. remove duplicates}
\For{$O_S\subseteq \texttt{Orb}_i$}
	\For{Every $F=F_1 \cup \cdots \cup F_{i-1} \cup O_S \cup F_{i+1} \cup \cdots \cup F_n$ with $F_j$ a $\lambda_j$-subset of the $j$-faces of $\Delta_n$}
	\State Compute a $\mathfrak S_{n+1}$ orbit representative $O_F$ of $F$.
	\EndFor
\EndFor
\State Define $\texttt{Orb}_\lambda$ to be the set of all representatives $O_F$ \Comment{i.e. remove duplicates}
\State\Return $\texttt{Orb}_\lambda$
\end{algorithmic}
\end{algorithm}

\subsection{Step 2: Determining which candidate orbits are bases} 
Once an exhaustive list of $\beta_n$-many candidate basis orbit representatives have been computed, we proceed to the task of determining which are, indeed, bases of $\mathcal M_n$. One algorithm for doing so is Algorithm \ref{alg:IsBasis}.
\begin{algorithm}[!htb]
\small
\caption{{\sf IsBasis}}\label{alg:IsBasis}
\algnotext{EndIf}
\algnotext{EndFor}
\algnotext{EndWhile}
\begin{flushleft}
\textbf{Input: }\\ \quad { $\varphi$: a polynomial parametrization  (over $\mathbb{Q}$) of an $r$-dimensional affine variety in $\mathbb{C}^N$ by $\mathbb{C}^r$}  \\
\quad \,{$B$: a subset of $[N]$ of size $r$} \\
\textbf{Output: }{$\texttt{true}$ if $B$ is a basis of $\mathcal M(X)$ and $\texttt{false}$ otherwise} \\
\end{flushleft}
\begin{algorithmic}[1]
\State Compute the submatrix $J_B$ of the Jacobian $J=\textrm{Jac}(\varphi)$ indexed by $B$
\State Compute the determinant $D_B$ of $J_B$ as a polynomial in $\mathbb{C}[x_1,\ldots,x_N]$
\If{$D_B$ is identically zero}{ \Return{\texttt{false}}}\Else{ \Return{\texttt{true}}}
\EndIf
\end{algorithmic}
\end{algorithm}

The main drawback of Algorithm \ref{alg:IsBasis} is that it becomes prohibitively expensive to compute on the order of  $\beta_4 = 48533$ many $10 \times 10$ determinants of matrices with entries in a polynomial ring of $e(4) = 10$ variables. Alternatively, we can evaluate the Jacobian at many points in $\mathbb{Q}^{10}$, giving the \textit{one-sided} Monte-Carlo algorithm,  Algorithm \ref{alg:IsBasisMonteCarlo}. The benefit of this approach is that determinants over $\mathbb{Q}$ are easier to compute than determinants over a polynomial ring.
\begin{algorithm}[!htb]
\small
\caption{{\sf IsBasis - One-sided Monte-Carlo (Randomized)}}\label{alg:IsBasisMonteCarlo}
\algnotext{EndIf}
\algnotext{EndFor}
\algnotext{EndWhile}
\begin{flushleft}
\textbf{Input: }\\ \quad { $\varphi$: a polynomial (over $\mathbb{Q}$) parametrization of an $r$-dimensional affine variety in $\mathbb{C}^N$ by $\mathbb{C}^r$}  \\
\quad \,{$B$: a subset of $[N]$ of size $r$} \\
\textbf{Output: }{$\texttt{true}$ if $B$ is a basis of $\mathcal M(X)$ and $\texttt{false}$ otherwise} \\
\end{flushleft}
\begin{algorithmic}[1]
\State Compute the submatrix $J_B$ of the Jacobian $J=\textrm{Jac}(\varphi)$ indexed by $B$
\State Evaluate $J_B$ at a random rational vector $q \in \mathbb{Q}^{r}$ to obtain $J_B(q)$.
\State Compute the determinant $D_B \in \mathbb{Q}$ of $J_B(q)$.
\If{$D_B=0$}{ \Return{\texttt{false}}}\Else{ \Return{\texttt{true}}}
\EndIf
\end{algorithmic}
\end{algorithm}

The probability of Algorithm \ref{alg:IsBasisMonteCarlo} outputting a false negative is no more than the probability that the rational number chosen lies on the \textit{non-matroidal locus} (see \cite{Rosen}) of the variety $X$. The non-matroidal locus is at most a hypersurface, so with a sufficient random generator for rational numbers, this probability is small. On the other hand, Algorithm \ref{alg:IsBasisMonteCarlo} never outputs a false positive, so it can always be used as a fast one-sided test. To increase the confidence of Algorithm \ref{alg:IsBasisMonteCarlo}, one may \textit{amplify} its likelihood of being correct by running it multiple times and outputting \texttt{false} if and only if \textit{every} output is \texttt{false}.

As a partner to Algorithm \ref{alg:IsBasisMonteCarlo}, we propose another one-sided algorithm in the other direction based on a popular solution-bound in computational algebraic geometry: the BKK bound.
\begin{proposition}[BKK bound \cite{Bernstein,Kushnirenko} ]
\label{prop:BKK}
Let $F = (f_1,\ldots,f_n) \subset \mathbb{C}[x_1,\ldots,x_n]$ be a square polynomial system with Newton polytopes $P = (P_1,\ldots,P_n)$. The number of isolated solutions to $F=0$ in the algebraic torus $(\mathbb{C}^\times)^n$ is bounded by the \textit{mixed volume} of $P$, known as the BKK bound of $F=0$. 
\end{proposition}
Proposition \ref{prop:BKK} may be adapted to bound the number of isolated solutions in affine space $\mathbb{C}^n$ (see \cite{LiWang}) which is the output of the \texttt{mixed\_volume} function in \texttt{HomotopyContinuation.jl} \cite{HCjl} with the option \texttt{only\_non\_zero} set to its default of \texttt{false}. We have the following one-sided algorithm.
\begin{algorithm}[!htb]
\small
\caption{{\sf IsBasis - One-sided BKK}}\label{alg:IsBasisBKK}
\algnotext{EndIf}
\algnotext{EndFor}
\algnotext{EndWhile}
\begin{flushleft}
\textbf{Input: }\\ \quad { $\varphi$: a polynomial parametrization (defined over $\mathbb{Q}$)  of an $r$-dimensional affine variety in $\mathbb{C}^N$ by $\mathbb{C}^r$}  \\
\quad \,{$B$: a subset of $[N]$ of size $r$} \\
\textbf{Output: }{$\texttt{false}$  only if $B$ is not a basis of $\mathcal M(X)$} \\
\end{flushleft}
\begin{algorithmic}[1]
\State Construct the square polynomial system $F_B = \{\varphi_i-x_i\}_{i \not\in B}\cup \{\varphi_i-b_i\}_{i \in B}$ where $b_i$ are taken to be generic coefficients.
\State Compute the (affine) mixed volume $\textrm{MV}_B$, bounding the number of isolated solutions of the system $F_B=0$ over $\mathbb{C}^N$ 
\If{$\textrm{MV}_B$ is zero}{ \Return{\texttt{false}}}\Else{ \Return{\texttt{nothing}}}
\EndIf
\end{algorithmic}
\end{algorithm}

\subsection{Step 3: Determining the degrees $d_B$} Given a basis $B$ of $\mathcal M_n$, the next step is to determine the degree $d_B$ of the branched cover $\pi_B:X \to \mathbb{C}_{B}^{e(n)}$. This amounts to solving a parametrized polynomial system. There are several ways to perform such a computation. We use techniques from numerical algebraic geometry and refer the reader to Section $3$ of \cite{MetricAG} for details.

\subsection{Step 4:  Determining the monodromy groups $G_{\pi_B}$}
We propose two algorithms for determining the monodromy groups $G_{\pi_B}$. The first is numerical.
\begin{algorithm}[!htb]
\small
\caption{{\sf Monodromy Group (Numerical and Randomized)}}\label{alg:NumericalMonodromyGroup}
\algnotext{EndIf}
\algnotext{EndFor}
\algnotext{EndWhile}
\begin{flushleft}
\textbf{Input: }\\
{$B$: a basis of the algebraic matroid of $X_n$}\\
{$F_B$: a parametrized polynomial system $F_B$ representing the branched cover $\pi_B:X_n \to \mathbb{C}_{B}^{|B|}$}   \\
\textbf{Options: }{$K$: the number of monodromy loops to take} \\
\textbf{Output: }{The monodromy group $G_{\pi_B}$} \\
\end{flushleft}
\begin{algorithmic}[1]
\State Sample several loops $\gamma_1,\ldots,\gamma_K$ in the parameter space $\mathbb{C}_{B}^{|B|}$
\State Compute a fibre $\pi_B^{-1}(\textbf{u})$ over some base parameter $\textbf{u}$ which is a regular value of $\pi_B$
\State Label the fibre $\pi_B^{-1}(\textbf{u})$ by the numbers $[d_B]$
\State Perform homotopy continuation over the loops $\gamma_1,\ldots,\gamma_K$ and compare the endpoints of each solution path to obtain the permutations $\sigma_{\gamma_1},\ldots,\sigma_{\gamma_K} \in \mathfrak S_{d_B}$
\State \Return{The subgroup of $\mathfrak S_{d_B}$ generated by $\sigma_{\gamma_1},\ldots,\sigma_{\gamma_K}$}
\end{algorithmic}
\end{algorithm}

The correctness of Algorithm \ref{alg:NumericalMonodromyGroup} is subject to two main concerns stemming from both its numerical and randomized nature. The concern pertaining to randomization is that the sampling procedure in \CHANGE{}{step} $\texttt{1}$ may not induce a reasonable distribution on $G_{\pi_B}$ and so the permutations computed may only produce a subgroup of the monodromy group. The methods in \cite{NumericalMonodromy} can mitigate this concern by computing a witness set for the discriminant of the enumerative problem and performing loops around the witness points. The numerical concerns appear in steps $\texttt{2}$ and $\texttt{4}$. Step $\texttt{2}$ is subject to the usual issues of any numerical solver, and step \texttt{4} is subject to the numerical errors involved in homotopy continuation. These numerical issues may be circumvented with a combination of \textit{certification} \cite{Certify}, degree bounds like Proposition \ref{prop:BKK}, and certified path-tracking (e.g. \cite{CertifiedPathTracking}). 

Assuming that there are no mistakes in the path-tracking in Algorithm \ref{alg:NumericalMonodromyGroup}, one is guaranteed to obtain a subgroup of $G_{\pi_B}$ as output. One may increase the confidence in the output of Algorithm \ref{alg:NumericalMonodromyGroup} by increasing the optional value $K$. For the other containment, we propose an algorithm (Algorithm \ref{alg:CoordinateSymmetryGroup}) for computing the coordinate symmetry group of $\pi_B$. 
\newpage

\begin{algorithm}[!htb]
\small
\caption{{\sf Coordinate Symmetry Group (Numerical and Randomized)}}\label{alg:CoordinateSymmetryGroup}
\algnotext{EndIf}
\algnotext{EndFor}
\algnotext{EndWhile}
\begin{flushleft}
\textbf{Input: }\\
{$B$: a basis of the algebraic matroid of $X_n$}\\
{$F_B$: a parametrized polynomial system $F_B$ representing the branched cover $\pi_B:X_n \to \mathbb{C}_{B}^{|B|}$}   \\
\textbf{Output: }{The coordinate symmetry group ${\widehat{G}}_{\pi_B}$} \\
\end{flushleft}
\begin{algorithmic}[1]
\State Compute a fibre $\pi_B^{-1}(\textbf{u})$ over some generic base parameter $\textbf{u}$ which is a regular value of $\pi_B$
\State Label the fibre $\pi_B^{-1}(\textbf{u})$ by the numbers $[d_B]$
\State For each coordinate $x_i$, partition the numbers $[d_B]$ based on whether the corresponding solutions agree in the $x_i$-th coordinate, label this partition $P_i$
\State Construct $\widehat{G}_{\pi_B} = \bigcap_{i=1}^n \mathfrak S_{P_i}$ as in Corollary \ref{cor:coordinatesymmetries}
\State \Return{$\widehat{G}_{\pi_B}$} 
\end{algorithmic}
\end{algorithm}

We remark that \texttt{Step 3} of Algorithm \ref{alg:CoordinateSymmetryGroup} relies on the ability to check if two complex numbers, each represented by a pair of floating point numbers, are equal. This is its numerical concern. From a randomized point-of-view, with probability zero a base parameter $\textbf{u}$ is chosen for which the points in its fibre have certain equal coordinates which occur from coincidence rather than from the structure of the problem. For example, the set of parameters $b,c$ for which the system $x^2+bx+c=0,y^2-b=0$ has a common coordinate in $x$ has measure zero.

\section{Results}
\label{sec:results}
We now collect our results pertaining to the algebraic matroids of the Heron varieties.

\subsection{The algebraic matroids of $X_2$, $X_3$, and $X_4$}
We begin with a result alluded to in Section \ref{sec:background} regarding {orbit representatives of candidate bases} for $\mathcal M_2$, $\mathcal M_3$, and $\mathcal M_4$. 
\begin{theorem}
Respectively, there are $(\beta_2,\beta_3,\beta_4) = (2, 35, 48533)$  orbits of subsets of positive dimensional faces of $(\Delta_2,\Delta_3,\Delta_4)$ of size $(e(2),e(3),e(4))=(3,6,10)$. 
\end{theorem}
\begin{proof}
After producing all unordered number partitions of $e(n)$ which are plausible $f$-vectors of $e(n)$-subsets of positive dimensional faces of $\Delta_n$, we apply Algorithm \ref{alg:OrbitsWithFvector} to achieve the result.
\end{proof}

We now describe the algebraic matroids of the first three Heron varieties. For $n=2$, this is essentially given by Example \ref{ex:X2Matroid}. We state it here for completeness.
\begin{theorem}[The algebraic matroid of $X_2$]
The algebraic matroid of the second Heron variety is the uniform matroid $U(3,4)$. The unique circuit is $\{12,13,23,123\}$.
\end{theorem}

\begin{figure}[!htpb]
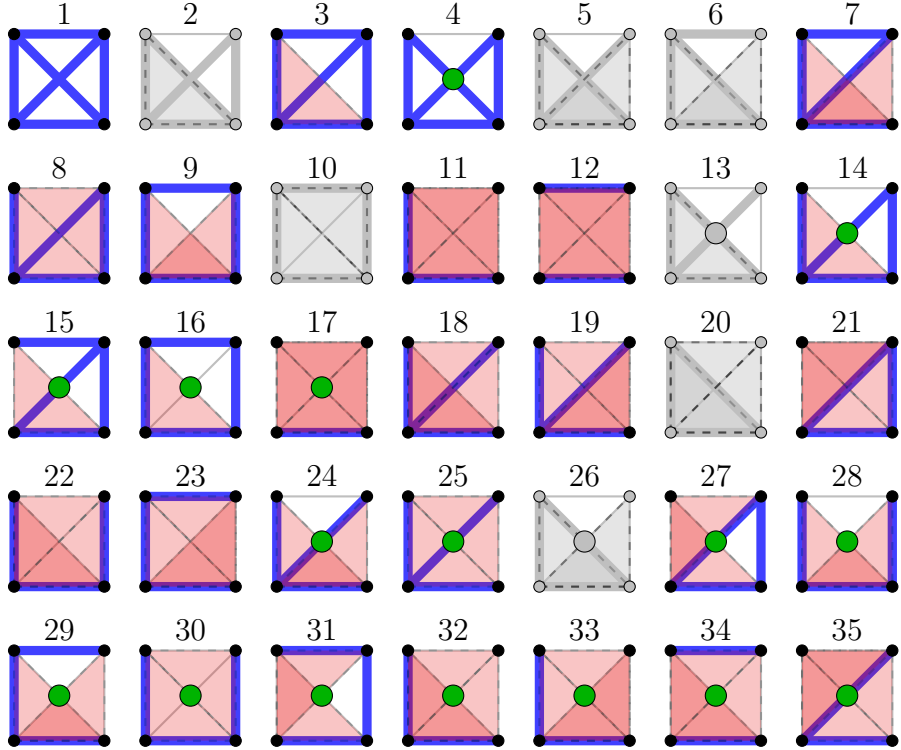

\centering

\OrbitsX{1/2,1/3,1/4,2/3,2/4,3/4}{}{} { $1$}
\OrbitsXNonbasis{1/2,1/3,1/4,2/3,2/4}{1/2/3}{} {$2$ }
\OrbitsX{1/2,1/3,1/4,2/4,3/4}{1/2/3}{} {$3$}
\OrbitsX{1/2,1/3,1/4,2/3,2/4}{}{1/2/3/4} {$4$}
\OrbitsXNonbasis{1/2,1/3,1/4,2/3}{1/2/3,1/2/4}{}  {$5$ }
\OrbitsXNonbasis{1/2,1/3,2/3,3/4}{1/2/3,1/2/4}{}  {$6$ }
\OrbitsX{1/2,1/3,1/4,3/4}{1/2/3,1/2/4}{} { $7$}

\vspace{0.2cm}

\OrbitsX{1/2,1/3,1/4,2/4}{1/2/3,2/3/4}{} {$8$}
\OrbitsX{1/2,1/3,2/4,3/4}{1/2/3,1/2/4}{} {$9$}
\OrbitsXNonbasis{1/2,1/3,2/4,3/4}{1/2/3,2/3/4}{} {$10$}
\OrbitsX{1/2,1/3}{1/2/3,1/2/4,1/3/4,2/3/4}{} {$11$}
\OrbitsX{1/2,3/4}{1/2/3,1/2/4,1/3/4,2/3/4}{} {$12$}
\OrbitsXNonbasis{1/2,1/3,1/4,2/3}{1/2/3}{1/2/3/4} {$13$}
\OrbitsX{1/2,1/3,1/4,2/4}{1/2/3}{1/2/3/4} { $14$}

\vspace{0.2cm}

\OrbitsX{1/2,1/4,2/4,3/4}{1/2/3}{1/2/3/4} {$15$}
\OrbitsX{1/2,1/3,2/4,3/4}{1/2/3}{1/2/3/4} {$16$}
\OrbitsX{1/2}{1/2/3,1/2/4,1/3/4,2/3/4}{1/2/3/4} {$17$}
\OrbitsX{1/2,1/3,1/4}{1/2/3,1/2/4,1/3/4}{} {$18$}
\OrbitsX{1/2,1/3,1/4}{1/2/3,1/2/4,2/3/4}{} {$19$}
\OrbitsXNonbasis{1/2,1/3,2/3}{1/2/3,1/2/4,1/3/4}{} {$20$}
\OrbitsX{1/2,1/4,2/4}{1/2/3,1/3/4,2/3/4}{} {$21$}

\vspace{0.2cm}

\OrbitsX{1/2,1/3,2/4}{1/2/3,1/2/4,1/3/4}{} {$22$}
\OrbitsX{1/2,1/3,3/4}{1/2/3,1/2/4,2/3/4}{} {$23$}
\OrbitsX{1/2,1/3,1/4}{1/2/3,1/2/4}{1/2/3/4} {$24$}
\OrbitsX{1/2,1/3,1/4}{1/2/3,2/3/4}{1/2/3/4} {$25$}
\OrbitsXNonbasis{1/2,1/3,2/3}{1/2/3,1/2/4}{1/2/3/4} {$26$}
\OrbitsX{1/2,1/4,2/4}{1/2/3,1/3/4}{1/2/3/4} {$27$}
\OrbitsX{1/2,1/3,2/4}{1/2/3,1/2/4}{1/2/3/4} {$28$}

\vspace{0.2cm}

\OrbitsX{1/2,1/3,3/4}{1/2/3,1/2/4}{1/2/3/4} {$29$}
\OrbitsX{1/2,1/3,2/4}{1/2/3,2/3/4}{1/2/3/4} {$30$}
\OrbitsX{1/2,2/4,3/4}{1/2/3,1/3/4}{1/2/3/4} {$31$}
\OrbitsX{1/2,1/3}{1/2/3,1/2/4,1/3/4}{1/2/3/4} {$32$}
\OrbitsX{1/2,1/3}{1/2/3,1/2/4,2/3/4}{1/2/3/4} {$33$}
\OrbitsX{1/2,3/4}{1/2/3,1/2/4,1/3/4}{1/2/3/4} {$34$}
\OrbitsX{1/2,1/4}{1/2/3,1/3/4,2/3/4}{1/2/3/4} {$35$}

\caption{The $35$ orbits of $6$-subsets of positive-dimensional faces of $\Delta_3$. Blue lines are edges, shaded red regions indicate triangular facets, and a green dot indicates the $3$-dimensional volume. The grey orbits represent non-bases of size $6$.} 
\label{fig:X3Picture}
\end{figure}

The matroid $\mathcal M_3$ is more interesting than $\mathcal M_2$. 
\begin{theorem}[The algebraic matroid of $X_3$]
\label{thm:X3Matroid}
The algebraic matroid $\mathcal M_3$ of the third Heron variety $X_3$ consists of $28$ orbits of bases and $7$ orbits of non-bases of size $e(3)=6$, under the action of $\mathfrak S_4$ on the vertices of $\Delta_3$. Each of these $35$ orbits is labeled and illustrated in Figure \ref{fig:X3Picture}. Those which are non-bases are in gray.

Up to the $\mathfrak S_4$-action, there are two orbits of circuits of size $6$ or smaller:
\[
C_1 = \{12,13,23,123\} \quad \quad \quad \text{ and } \quad \quad \quad C_2 =\{12,13,24,34,123,234\}.
\]
\end{theorem}
\begin{proof}
We apply Algorithm \ref{alg:IsBasis} to each of the $35$  orbit representatives of candidate bases of $\mathcal M_3$ (obtained via Algorithm \ref{alg:OrbitsWithFvector}) to determine which are bases and obtain the result.
\end{proof}
Six of the seven orbits of non-bases of $\mathcal M_3$ of size $6$ are easily seen to consist of non-bases: each representative contains a  $\mathfrak S_4$-representative of the circuit $\{12,13,23,123\}$ of $\mathcal M_2$.  In this way, dependent subsets of $\mathcal M_n$ induce dependent subsets of $\mathcal M_{n+k}$ in the following sense.
\begin{lemma}
\label{lem:circuitinduction}
If $D$ is dependent in $\mathcal M_n$, then for any $k\CHANGE{}{\in \mathbb{N}}$, $D$ is dependent in $\mathcal M_{n+k}$ as well. 
\end{lemma}
\begin{proof}
Consider the forgetful map $f_n:X_{n+k} \to X_{n}$ from the $(n+k)$-th Heron variety to the $n$-th Heron variety which simply forgets all coordinates involving a vertex labeled $n+2$ or higher.
This map is clearly dominant as every $n$-simplex can be realized as an $n$-face of an $(n+k)$-simplex. Fix a dependent set $D$ of $\mathcal M_n$, that is, a subset $D$ of the ground set of $\mathcal M_n$ such that $\pi_D:X_n \to \mathbb{C}_D^{|D|}$ is not dominant. Since $\pi'_D:X_{n+k} \to \mathbb{C}_D^{|D|}$ factors as $\pi'_D = \pi_D \circ f_n$, we conclude that $\pi'_D$ is also not dominant.
\end{proof}
\begin{remark}
\label{remark:circuitinductionfails}
The circuit $B_{10} = \{12,13,24,34,123,234\}$ of $\mathcal M_3$ is not implied by Lemma \ref{lem:circuitinduction}, however, a rigidity argument shows that $B_{10}$ is indeed dependent. Suppose that $\Delta_3(\textbf{p})$ is a simplex and all faces other than those contained in the two triangular faces $123$ and $234$ are deleted. These two triangles share one edge, and there is a rigid motion where the common edge functions as a hinge. Along this motion, the edge length $v_{14}$ changes, and so the six volumes of $B_{10}$ cannot determine the remaining edge lengths. One may also realize $B_{10}$ as a non-basis via the \textit{circuit exchange axiom} on the circuits $\{12,13,23,123\}$ and $\{23,24,34,234\}$.
\end{remark}

As shown in Remark \ref{remark:circuitinductionfails}, Lemma \ref{lem:circuitinduction} does not identify every dependent set in $\mathcal M_3$, however, our next result states that Algorithm \ref{alg:IsBasisBKK} does.
\begin{proposition}
A subset $B$ of size $6$ of the ground set of $\mathcal M_3$ is a non-basis if and only if the mixed volume of the polynomial system describing a fibre of $\pi_B$ is zero. In particular, Algorithm \ref{alg:IsBasisBKK} characterizes the bases of $\mathcal M_3$. 
\end{proposition}
\begin{proof}
We constructed each of the polynomial systems corresponding to the $35$ orbit representatives and computed their mixed volume. We observed that an orbit representative had non-zero mixed volume if and only if it was one of the bases of Theorem \ref{thm:X3Matroid}.
\end{proof}

\begin{remark}
We remark that given any irreducible complex variety $X \subseteq \mathbb{C}^N$ of dimension $r$ cut out by $N-r$ polynomials $f_1,\ldots,f_{N-r}$, we may try to compute non-bases of $\mathcal M(X)$ using Algorithm \ref{alg:IsBasisBKK}. Those polynomials $f_1,\ldots,f_{N-r}$ for which Algorithm \ref{alg:IsBasisBKK} characterizes the bases of their vanishing locus are special and deserve further research: such polynomial systems have ``generic algebraic matroids'' in terms of their sparse structure.
\end{remark}

\begin{theorem}[The algebraic matroid of $X_4$]
\label{thm:matroidX4}
Of the $\beta_4=48533$ orbit representatives of $10$-subsets of positive-dimensional faces of $\Delta_4$, $35887$ of them are bases and $12646$ are not.
\end{theorem}
\begin{proof}
First, the one-sided test of Algorithm \ref{alg:IsBasisMonteCarlo} is applied to the $48533$ candidates to obtain some bases. In fact, this successfully finds all bases. Algorithm \ref{alg:IsBasisBKK} is applied to the remaining candidates to certify non-bases. After these two fairly quick procedures, $155$ candidates remain to be classified. Applying Algorithm \ref{alg:IsBasis} to these candidates confirms that all are non-bases and comprises the most demanding part of the procedure, but is still completed in under twenty minutes. 
\end{proof}

Unlike the case where $n=3$, neither Lemma \ref{lem:circuitinduction} nor  Algorithm \ref{alg:IsBasisBKK} characterizes the non-bases of $\mathcal M_4$. We quantify the extent to which this is the case below.

\begin{corollary}
Of the $12646$ non-basis orbits of $\mathcal M_4$, Lemma \ref{lem:circuitinduction} fails to identify $1166$ of them and Algorithm \ref{alg:IsBasisBKK} fails to identify a \textit{strictly} smaller set of $155$ of them.  
\end{corollary}

\begin{example}
Consider the subsets
\[S = \{12,14,15,24,25,123,145,1234,1235,12345\}\]
\[
 S' = \{12,14,15,24,25,45,123,1234,1235,12345\}\]
 of positive-dimensional faces of $\Delta_4$ illustrated in Figure \ref{fig:sprime}.
 
\begin{figure}[!hptb]
\includegraphics[scale=0.5]{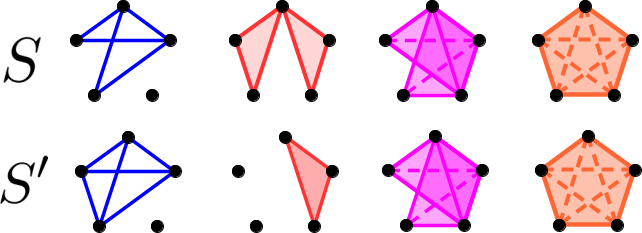}
\caption{Illustrations of two subsets of positive dimensional faces of $\Delta_4$. The edges are blue, the triangles are red, the tetrahedra are pink, and the $4$-simplex is orange.}
\label{fig:sprime}
\end{figure} 

  The subset $S$  
 does not contain a circuit of $\mathcal M_2$ or $\mathcal M_3$ realized within a facet of $\Delta_4$. However, the mixed volume of the polynomial system corresponding to $S$ is zero, and hence, it is a non-basis by Algorithm \ref{alg:IsBasisBKK}. The subset $S'$ on the other hand
 has a non-zero mixed volume of $16$ (and does not contain any circuit of $\mathcal M_2$ or $\mathcal M_3$). Nonetheless, it too is a non-basis of size ten, as confirmed by the computational justification of Theorem \ref{thm:matroidX4}.
\end{example}

\subsection{The decorations $d_B$ and $G_{\pi_B}$ of the bases of $\mathcal M_2$ and $\mathcal M_3$}
Solving a single instance of the parametrized polynomial system corresponding to a basis of $\mathcal M_n$ is simple using out-of-the box homotopy continuation software like \texttt{HomotopyContinuation.jl} \cite{HCjl}. Such a computation gives the base degree $d_B$. Beyond such computations, we also compute the monodromy groups of the branched covers corresponding to these bases. 
\begin{theorem}
\label{thm:X2Galois}
The monodromy groups associated to the two orbits $\{12,13,23\}$ and $\{12,13,123\}$ of bases of $\mathcal M_2$ are $\mathfrak S_1$ and $\mathfrak S_2$, respectively. 
\end{theorem}
\begin{proof}
The proof of Theorem \ref{thm:X2Galois} is trivial: Example \ref{ex:X2Matroid} already determined the base degrees of the bases of $\mathcal M_2$ and the monodromy groups must be transitive since $X_2$ is irreducible. 
\end{proof}

The following theorem provides an upper bound on the monodromy/Galois groups of the branched covers associated to bases of $X_3$.

\begin{figure}[!htpb]
\includegraphics[scale=0.4]{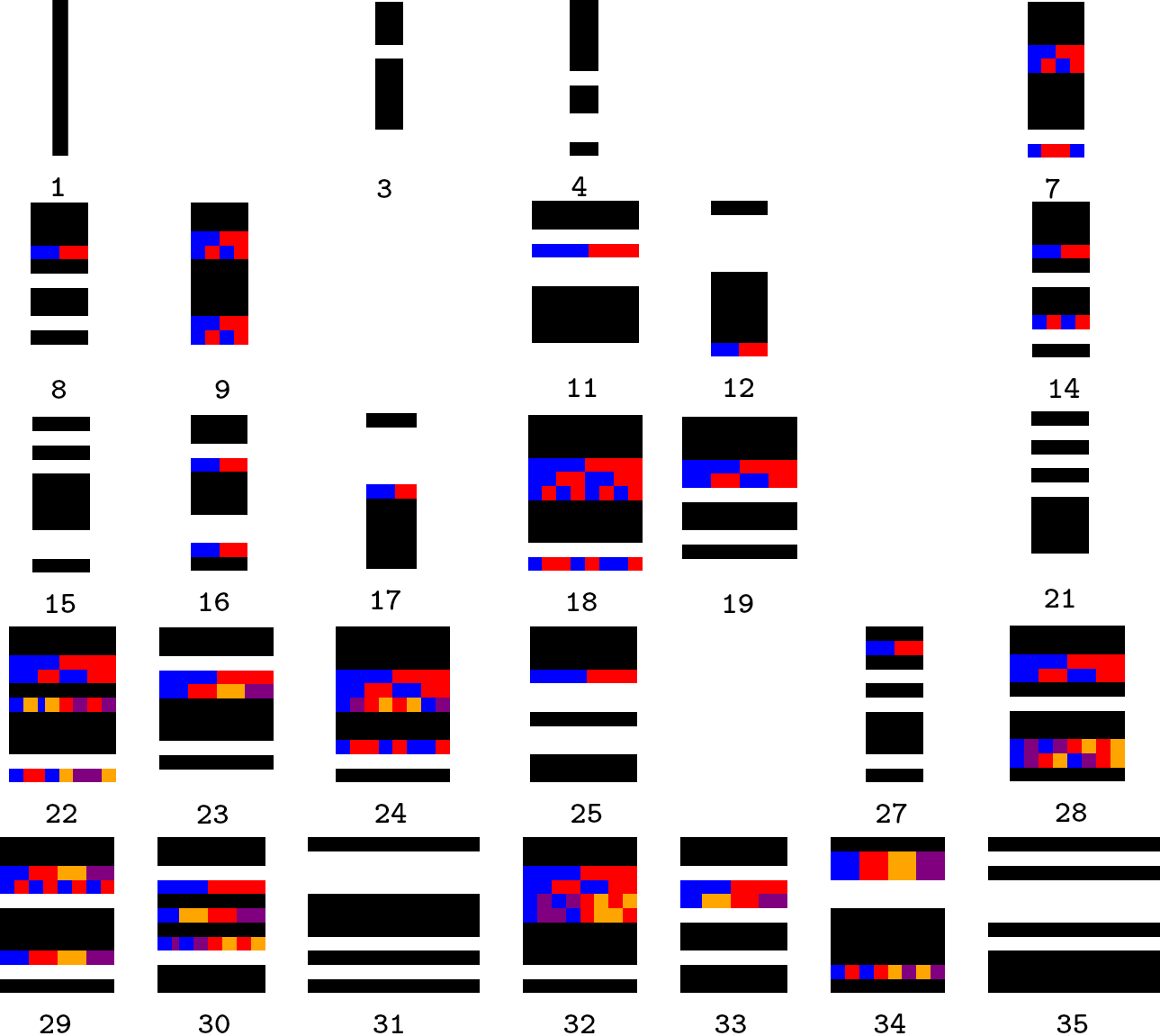}
\caption{A coordinate symmetry diagram for each basis  $B$ of $\mathcal M_3$ labeled according to Figure \ref{fig:X3Picture}. The $d_B$ columns correspond to points in the fibre of $\pi_B$ and the rows are labeled, in order, as $\{12,13,14,23,24,34,123,124,134,234,1234\}$. Black pixels indicate that every coordinate is the same and white pixels indicate that no coordinate is the same. }
\label{fig:symmetrydiagrams}
\end{figure}

\begin{theorem}[*\footnote{Asterisks are popularly employed in the field of numerical algebraic geometry to indicate that a result relies on the absence of any numerical error, or on a probability $1$ condition holding. The potential sources of errors for the results in this paper are written next to the name of each algorithm used and discussed within the body of the text.}]
\label{thm:CoordinateSymmetryX3}
The coordinate \CHANGE{monodromy}{symmetry} groups of the branched covers corresponding to bases of $\mathcal M_3$ are given by the coordinate symmetry diagrams in Figure \ref{fig:symmetrydiagrams}. See Table \ref{tab:GaloisGroups} for the descriptions of the isomorphism classes of these groups.
\end{theorem}
\begin{proof}
This result is obtained in two steps. The first is the solving of a generic fibre of $\pi_B$. After this step, we determine numerically which coordinates of the points in this fibre are equal. This essentially gives the diagrams in Figure \ref{fig:symmetrydiagrams}. To obtain the group descriptions in the fourth row of Table \ref{tab:GaloisGroups}, for each basis, we intersect the wreath products which preserve the blocks of each row of the coordinate symmetry diagram associated to that basis.  This is the content of Algorithm \ref{alg:CoordinateSymmetryGroup}. After performing this intersection, we use the \texttt{GAP} command \texttt{describe} to obtain a description of the isomorphism classes of the groups computed.
\end{proof}

\begin{table}[!htpb]
{\small{
\begin{tabular}{|c||c|c|c|c|c|c|c|c|c|c|c|c|c|c|c|}\hline
$i$ & $ \textbf{1} $ & $  \textbf{3} $ & $  \textbf{4} $ & $  \textbf{7} $ & $  \textbf{8} $ & $  \textbf{9} $ & $  \textbf{11} $ & $  \textbf{12} $& $  \textbf{14} $ & $  \textbf{15} $ & $  \textbf{16} $ & $  \textbf{17} $ & $  \textbf{18} $ & $  \textbf{19} $ & $  \textbf{21} $  \\ \hline 
$d_B$ & $1$ & $2$ & $ 2$ & $ 4$ & $ 4$ & $ 4$ & $ 8$ & $ 4$ & $ 4$ & $ 4$ & $ 4$ & $ 4$ & $ 8$ & $ 8$ & $ 4$\\ \hline
$G_{\pi_B}$& $\{e\}$ & $\ZZZ$ & $\ZZZ$ & $V$  & $\DDD$ & $V$ & $\mathfrak S_4 \wr \ZZZ$  &$\DDD$& $V$ &\alert{$\DDD$}& $\DDD$ & $\DDD$ & $\ZZZ^3$ &$\ZZZ \wr V$ & $\mathfrak S_4$\\ \hline
$\widehat{G}_{\pi B}$ & $\{e\}$ & $\ZZZ$ & $\ZZZ$ & $\Klein$ & $\DDD$ & $\Klein$ & $\mathfrak S_4 \wr \ZZZ$ & $\DDD$ & $\Klein$ & $\mathfrak S_4$ & $\DDD$ & $\DDD$ & $\ZZZ^3 $ & $\ZZZ \wr V $ & $\mathfrak S_4$ \\ \hline \hline 
\end{tabular}}}
\begin{tabular}{|c||c|c|c|c|c|c|c|c|}\hline
& $  \textbf{22}$ & $  \textbf{23} $ & $  \textbf{24} $ & $  \textbf{25} $ & $  \textbf{27} $ & $  \textbf{28} $ & $  \textbf{29} $ & $  \textbf{30} $  \\ \hline 
$d_B$ & $ 8$ &  $ 8$ & $ 8$ & $ 8$ & $ 4$ & $ 8$ & $ 8$ & $ 8$  \\ \hline 
$G_{\pi_B}$ &$\ZZZ \times \DDD$  & $\ZZZ \wr D_8$ & $(\ZZZ)^3$ & $\alert{\DDD \wr \ZZZ}$ & $\alert{V}$ & $(\ZZZ)^3$ &$\alert{\ZZZ \times \DDD}$ &$V \wr \ZZZ$ \\ \hline 
$\widehat{G}_{\pi B}$ &  $\ZZZ \times \DDD$ & $\ZZZ \wr D_8$  & $(\ZZZ)^3$ & $\mathfrak S_4 \wr \ZZZ$ & $\DDD$ & $(\ZZZ)^3$ & $\ZZZ \times \mathfrak S_4$ & $V \wr \ZZZ$ \\ \hline \hline 
\end{tabular}
\begin{tabular}{|c||c|c|c|c|c|}\hline
$i$ &$  \textbf{31} $& $  \textbf{32} $ & $  \textbf{33} $ & $  \textbf{34} $ & $  \textbf{35} $ \\ \hline
$d_B$ &$ 12$ & $ 8$ & $ 8$ & $ 8$ & $ 12$ \\ \hline 
$G_{\pi_B}$ &$\mydef{\mathfrak S_{12}}$ &  $\alert{(\ZZZ)^3}$ &\alert{$\ZZZ \times \DDD$}&$V \wr \ZZZ$ & $\mydef{\mathfrak S_{12}}$\\ \hline
$\widehat{G}_{\pi B}$ & $\mydef{\mathfrak S_{12}}$ & $\ZZZ \times \DDD$ & $\ZZZ \wr D_8$ & $V \wr \ZZZ$ & $\mydef{\mathfrak S_{12}}$  \\ \hline
\end{tabular}
\caption{The index, degree, numerically computed monodromy group, and coordinate symmetry group corresponding to each basis of $\mathcal M_3$, indexed as in Figure \ref{fig:X3Picture}. The numerically computed monodromy groups which differ from the coordinate symmetry groups are highlighted in red. Those groups which are not solvable are highlighted in blue. Here $\mydef{V}$ is the group $\ZZZ \times \ZZZ$.}
\label{tab:GaloisGroups}
\end{table}

\begin{theorem}[*]\label{thm:GaloisX3}  Table \ref{tab:GaloisGroups} lists the numerically computed degrees and isomorphism classes of monodromy groups of the branched covers corresponding to bases of $\mathcal M_3$ as indexed by Figure \ref{fig:X3Picture}. All permutation groups in Table \ref{tab:GaloisGroups} which are not the full symmetric group are imprimitive.
\end{theorem}
\begin{proof}
This theorem is obtained by applying Algorithm \ref{alg:NumericalMonodromyGroup} to a representative of each basis of $\mathcal M_3$. 
\end{proof}

One major assumption which is often made when computing monodromy groups numerically is that the sampling procedure for loops in the parameter space induces a somewhat uniform sample on the monodromy group. Evidence for the validity of this assumption is sparse, and so we include an example illustrating the distribution we encountered for $\pi_{B_{33}}$

\begin{example}
\label{ex:poormonodromydistribution}
Our sampling procedure for loops in $P = \mathbb{C}_{B_{33}}^{6}$ is induced by a sampling procedure for points in that same parameter space. After computing a fixed generic fibre of $\pi_{B_{33}}$ we sample two additional points in $P$ to make a triangle $\gamma$. Sampling points in $P$ is done coordinate-wise (real and imaginary separately) under the normal distribution of mean zero and standard deviation one. Optionally, we have a parameter $\texttt{Radius}$ which is a multiplier for the resulting points. Figure \ref{fig:GaloisDistribution} shows, under various radii, the permutations in $G_{\pi_{B_{33}}}$ which were found. 
\begin{figure}[!htpb]
\includegraphics[scale=0.35]{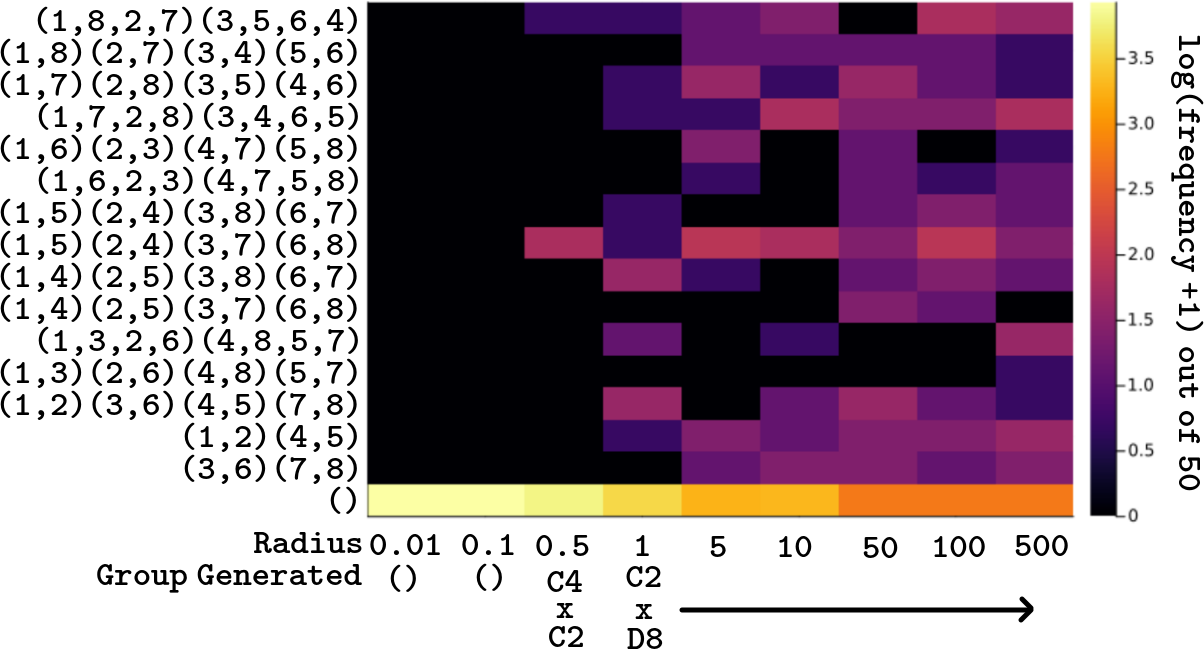}
\caption{A heatmap whose $(r,\sigma)$-th entry represents the frequency of observing the permutation $\sigma$ under the loop sampling procedure described with radius $r$. For each \CHANGE{radii}{radius}, we sampled $50$ loops. For visual purposes the heatmap value is the logarithm of $1$ plus the frequency.}
\label{fig:GaloisDistribution}
\end{figure}
\end{example}

Even though there is a gap between some of the  numerically computed monodromy groups in Table \ref{tab:GaloisGroups} and the coordinate symmetry groups, all groups in that table other than $\mathfrak S_{12}$ are solvable, and those which are not the full symmetric group are imprimitive. For those projections with solvable Galois groups, every coordinate in a fibre can be solved for by radicals in terms of the base coordinates. For the bases $B_{31}$ and $B_{35}$, not all coordinates can be solved for by radicals in the base coordinates. Moreover, since the coordinate symmetry groups of $B_{31}$ and $B_{35}$ are trivial, \textit{none} of the coordinates of fibres of $\pi_{B_{31}}$ or $\pi_{B_{35}}$ can be solved for by radicals in the base coordinates. We summarize this in the following corollary.

\begin{figure}[!b]
\includegraphics[scale=0.5]{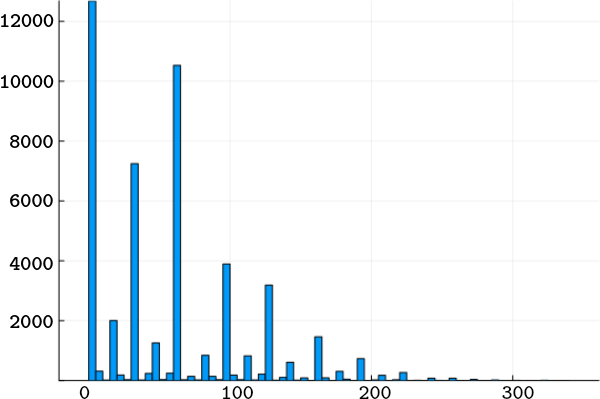}
\caption{A histogram of the base degrees associated to the $48533$ orbits of $10$-subsets of $\mathcal M_4$. A base degree of $0$ encodes that the corresponding orbit is not actually a basis. The maximum base degree is $336$.}
\label{fig:X4Degrees}
\end{figure}

\begin{corollary}[*]
\label{cor:solvability}
Fix a basis $B \in \mathcal M_3$ other than 
\[
B_{31} = \{12,13,123,124,134,1234\} \quad \quad \quad \quad \text{ or } \quad \quad \quad \quad B_{35} = \{12,14,123,134,234,1234\}.
\]
All volumes of a $3$-simplex are solvable in radicals as functions of the volumes indexed by $B$. Conversely, there is no such formula for the remaining volumes of a $3$-simplex for $B_{31}$ and $B_{35}$.
\end{corollary}
All of the methods described in Section \ref{sec:pipeline} may be applied to compute the coordinate symmetry groups and monodromy groups associated to bases of $\mathcal M_4$ as well. We performed these calculations for the base degrees, but not the Galois groups. A single example of a Galois group associated to a basis of $\mathcal M_4$ is given in Section \ref{secsec:physics}. A histogram of the base degrees of the $48533$ orbits of $10$-subsets of the ground set of $\mathcal M_4$ \CHANGE{are}{is} given in Figure \ref{fig:X4Degrees}. These numbers range from $0$ (indicating a non-basis) to $336$.

\section{Experiments}
\label{sec:experiments}
In this section, we analyze the polynomial systems arising from bases of $\mathcal M_3$ experimentally. For each basis $B$ of $\mathcal M_3$ we sample values for each volume $v_{s}$ associated to $s \in B$ uniformly at random within the unit interval $[0,1]$ and compute the fibre of $\pi_B$ over $\textbf{x}_B = \{v_{s}^2\}_{s \in B}$. We count the number of points in each fibre contained in the images of $\mathcal E_3$, $\mathbb{R}_{\CHANGE{\geq}{>0}}^3$, and $\mathbb{R}^3$. We do this for $100,000$ fibres and collect the data in Table \ref{tab:experimentsX3}.

\begin{table}
{\small{
\begin{tabular}{|c||c|c|c|} \hline 
$(i,d_{B_i})$ & $(1,\mydef{1})$ & $(3,\mydef{2})$ & $(4,\mydef{2})$ \\ \hline  
 \begin{tabular}{r}  
Data \\ \hline 
Real\\
Positive \\
Realizable \\  
\end{tabular}
& 
 \begin{tabular}{c|c}  
  Observed & Expected \\  
$\{\mydef{1}\}$ & 1.000 \\
$\{\mydef{1}\}$ & 1.000 \\
$\{0,\mydef{1}\}$ & 0.102  \\  
 \end{tabular}
& 
 \begin{tabular}{c|c}  
  Observed & Expected \\ \hline 
$\{0,\mydef{2}\}$ & 0.248\\
 $\{0,\mydef{2}\}$ & 0.248 \\
 $\{0,1,\mydef{2}\}$ & 0.031 \\  
 \end{tabular}
& 
 \begin{tabular}{c|c}  
  Observed & Expected \\ \hline 
$\{0,\mydef{2}\}$ &  0.048 \\
$\{0,1,\mydef{2}\}$ & 0.029 \\
$\{0,\mydef{2}\}$ & 0.009 \\  
 \end{tabular}\\ \hline \hline 
 
$(i,d_{B_i})$  & $(7,\mydef{4})$ & $(8,{\color{violet}{4}})$ & $(9,\mydef{4})$ \\ \hline  
 \begin{tabular}{r}  
Data \\ \hline 
Real\\
Positive \\
Realizable \\  
\end{tabular}
& 
 \begin{tabular}{c|c}  
  Observed & Expected \\ \hline 
$\{0,\mydef{4}\}$ &0.083  \\
$\{0,\mydef{4}\}$ & 0.083 \\
$\{0,2,\mydef{4}\}$ & 0.014 \\  
 \end{tabular}
& 
 \begin{tabular}{c|c}  
  Observed & Expected \\ \hline 
$\{0,2,\mydef{4}\}$ &0.108  \\
$\{0,2,\mydef{4}\}$ &0.108  \\
$\{0,\ldots,3\}$ & 0.014  \\  
 \end{tabular}
& 
 \begin{tabular}{c|c}  
  Observed & Expected \\ \hline 
$\{0,\mydef{4}\}$ & 0.082 \\
$\{0,\mydef{4}\}$ & 0.082 \\
$\{0,\ldots,\mydef{4}\}$ & 0.011 \\  
 \end{tabular}\\ \hline \hline 
 
$(i,d_{B_i})$  & $(11,{\color{violet}{8}})$ & $(12,\mydef{4})$ & $(14,\mydef{4})$ \\ \hline  
 \begin{tabular}{r}  
Data \\ \hline 
Real\\
Positive \\
Realizable \\  
\end{tabular}
& 
 \begin{tabular}{c|c}  
  Observed & Expected \\ \hline 
$\{0,4,6,\mydef{8}\}$ & 0.515 \\
$\{0,4,6,\mydef{8}\}$ & 0.515 \\
$\{0,2,4,6\}$ & 0.241  \\  
 \end{tabular}
& 
 \begin{tabular}{c|c}  
  Observed & Expected \\ \hline 
$\{0,2,\mydef{4}\}$ & 0.876 \\
$\{0,2,\mydef{4}\}$ & 0.876 \\
$\{0,2,\mydef{4}\}$ &0.059  \\  
 \end{tabular}
& 
 \begin{tabular}{c|c}  
  Observed & Expected \\ \hline 
$\{0,\mydef{4}\}$ & 0.005 \\
$\{0,\mydef{4}\}$ & 0.005 \\
$\{0,\mydef{4}\}$ & 0.005 \\  
 \end{tabular}\\ \hline \hline 
 
$(i,d_{B_i})$  & $(15,\mydef{4})$ & $(16,\mydef{4})$ & $(17,\mydef{4})$ \\ \hline  
 \begin{tabular}{r}  
Data \\ \hline 
Real\\
Positive \\
Realizable \\  
\end{tabular}
& 
 \begin{tabular}{c|c}  
  Observed & Expected \\ \hline 
$\{0,2,\mydef{4}\}$ &  0.005 \\
$\{0,2,\mydef{4}\}$ & 0.005 \\
$\{0,2,\mydef{4}\}$ & 0.005 \\  
 \end{tabular}
& 
 \begin{tabular}{c|c}  
  Observed & Expected \\ \hline 
$\{0,2,\mydef{4}\}$ & 0.003 \\
$\{0,2,\mydef{4}\}$ & 0.003 \\
$\{0,2,\mydef{4}\}$ & 0.003 \\  
 \end{tabular}
& 
 \begin{tabular}{c|c}  
  Observed & Expected \\ \hline 
$\{0,2,\mydef{4}\}$ & 0.081  \\
$\{0,2,\mydef{4}\}$ & 0.081  \\
$\{0,2,\mydef{4}\}$ &  0.081 \\  
 \end{tabular}\\ \hline \hline 
 
$(i,d_{B_i})$  & $(18,\mydef{8})$ & $(19,{\color{violet}{8}})$ & $(21,\mydef{4})$ \\ \hline  
 \begin{tabular}{r}  
Data \\ \hline 
Real\\
Positive \\
Realizable \\  
\end{tabular}
& 
 \begin{tabular}{c|c}  
  Observed & Expected \\ \hline 
$\{0,\mydef{8}\}$ & 0.037  \\
$\{0,\mydef{8}\}$ & 0.037 \\
$\{0,4,\mydef{8}\}$ & 0.008 \\  
 \end{tabular}
& 
 \begin{tabular}{c|c}  
  Observed & Expected \\ \hline 
$\{0,2,4,6,\mydef{8}\}$ & 0.064 \\
$\{0,2,4,6,\mydef{8}\}$ & 0.064 \\
$\{0,\ldots,5\}$ & 0.009 \\  
 \end{tabular}
& 
 \begin{tabular}{c|c}  
  Observed & Expected \\ \hline 
$\{0,2,\mydef{4}\}$ & 1.076 \\
$\{0,2,\mydef{4}\}$ & 1.076 \\
$\{0,\ldots,\mydef{4}\}$ & 0.472 \\  
 \end{tabular}\\ \hline \hline 
 
$(i,d_{B_i})$  & $(22,{\color{violet}{8}})$ & $(23,{\color{violet}{8}})$ & $(24,\mydef{8})$ \\ \hline  
 \begin{tabular}{r}  
Data \\ \hline 
Real\\
Positive \\
Realizable \\  
\end{tabular}
& 
 \begin{tabular}{c|c}  
  Observed & Expected \\ \hline 
$\{0,4,\mydef{8}\}$ & 0.049  \\
$\{0,4,\mydef{8}\}$ & 0.049 \\
$\{0,2,4,6\}$ & 0.009 \\  
 \end{tabular}
& 
 \begin{tabular}{c|c} \hline 
  Observed & Expected \\ \hline 
$\{0,2,4,6,\mydef{8}\}$ & 0.098 \\
$\{0,2,4,6,\mydef{8}\}$ & 0.098 \\
$\{0,1,2,3,4\}$ & 0.014 \\  
 \end{tabular}
& 
 \begin{tabular}{c|c} \hline 
  Observed & Expected \\ \hline 
$\{0,\mydef{8}\}$ & 0.003 \\
$\{0,\mydef{8}\}$ & 0.003 \\
$\{0,\mydef{8}\}$ & 0.003 \\  
 \end{tabular}\\ \hline \hline 
 
$(i,d_{B_i})$  & $(25,\mydef{8})$ & $(27,\mydef{4})$ & $(28,\mydef{8})$ \\ \hline  
 \begin{tabular}{r}  
Data \\ \hline 
Real\\
Positive \\
Realizable \\  
\end{tabular}
& 
 \begin{tabular}{c|c}  
  Observed & Expected \\ \hline 
$\{0,2,4,6,\mydef{8}\}$ & 0.004  \\
$\{0,2,4,6,\mydef{8}\}$ & 0.004  \\
$\{0,2,4,6,\mydef{8}\}$ & 0.004  \\  
 \end{tabular}
& 
 \begin{tabular}{c|c}  
  Observed & Expected \\ \hline 
$\{0,\mydef{4}\}$ & 0.078 \\
$\{0,\mydef{4}\}$ & 0.078  \\
$\{0,\mydef{4}\}$ & 0.078 \\  
 \end{tabular}
& 
 \begin{tabular}{c|c}  
  Observed & Expected \\ \hline 
$\{0,\mydef{8}\}$ & 0.002 \\
$\{0,\mydef{8}\}$ & 0.002 \\
$\{0,\mydef{8}\}$ & 0.002 \\
 \end{tabular}\\ \hline \hline 
 
$(i,d_{B_i})$  & $(29,\mydef{8})$ & $(30,\mydef{8})$ & $(31,\alert{12})$ \\ \hline  
 \begin{tabular}{r}  
Data \\ \hline 
Real\\
Positive \\
Realizable \\  
\end{tabular}
& 
 \begin{tabular}{c|c}  
  Observed & Expected \\ \hline 
$\{0,4,\mydef{8}\}$ & 0.002 \\
$\{0,4,\mydef{8}\}$ & 0.002 \\
$\{0,4,\mydef{8}\}$ & 0.002 \\
 \end{tabular}
& 
 \begin{tabular}{c|c}  
  Observed & Expected \\ \hline 
$\{0,4,\mydef{8}\}$ & 0.003 \\
$\{0,4,\mydef{8}\}$ & 0.003 \\
$\{0,4,\mydef{8}\}$ & 0.003 \\
 \end{tabular}
& 
 \begin{tabular}{c|c}  
  Observed & Expected \\ \hline 
$\{0,2,4,6,8\}$ & 0.006  \\
$\{0,2,4,6,8\}$ & 0.006  \\
$\{0,2,4,6,8\}$ & 0.006  \\
 \end{tabular}\\ \hline \hline 
 
$(i,d_{B_i})$ & $(32,\mydef{8})$ & $(33,\mydef{8})$ & $(34,\mydef{8})$ \\ \hline  
 \begin{tabular}{r}  
Data \\ \hline 
Real\\
Positive \\
Realizable \\  
\end{tabular}
& 
 \begin{tabular}{c|c}  
  Observed & Expected \\ \hline 
$\{0,\mydef{8}\}$ & 0.035 \\
$\{0,\mydef{8}\}$ & 0.035 \\
$\{0,\mydef{8}\}$ & 0.035 \\
 \end{tabular}
& 
 \begin{tabular}{c|c}  
  Observed & Expected \\ \hline 
$\{0,4,\mydef{8}\}$ & 0.0304 \\
$\{0,4,\mydef{8}\}$ & 0.0304 \\
$\{0,4,\mydef{8}\}$ & 0.0304 \\
 \end{tabular}
& 
 \begin{tabular}{c|c}  
  Observed & Expected \\ \hline 
$\{0,4,\mydef{8}\}$ & 0.039 \\
$\{0,4,\mydef{8}\}$ & 0.039 \\
$\{0,4,\mydef{8}\}$ & 0.039 \\
 \end{tabular}\\ \hline \hline 
 
$(i,d_{B_i})$ & $(35,\mydef{12})$ &  & \\ \hline  
 \begin{tabular}{r}  
Data \\ \hline 
Real\\
Positive \\
Realizable \\  
\end{tabular}
& 
 \begin{tabular}{c|c}  
  Observed & Expected \\ \hline 
$\{0,2,4,6,8,10,\mydef{12}\}$ & 0.079 \\
$\{0,2,4,6,8,10,\mydef{12}\}$ & 0.079 \\
$\{0,2,4,6,8,10,\mydef{12}\}$ & 0.079 \\
 \end{tabular}
& 

& 
\\ \hline \hline

\end{tabular}}}
\caption{A table which summarizes experimental data. We computed $100,000$ fibres of $\pi_{B_i}$ for each $i\in\{1,\ldots,35\}$ which indexes a basis of $\mathcal M_3$ and tabulated our observations. The degree $d_B$ is blue if a fibre which is completely realizable was found, purple if a fibre which is completely positive was found, and red if no fibre was found that is completely real. }
\label{tab:experimentsX3}
\end{table}

Equipped with this experimental data, we proceed to make several observations. Keep in mind that the parameters we chose were special in a semi-algebraic sense: for any basis $B$, we computed fibres exclusively over points in $[0,1]^6 \subset \mathbb{C}_{B}^{|B|}$. 

\subsection{Congruence modulo four:} Several of the enumerative problems $\pi_B:X_3 \to \mathbb{C}_{B}^{|B|}$ for bases $B$ of $\mathcal M_3$ appear to have real solution counts which are equivalent modulo four. Specifically, these bases are 
\begin{gather*}
\label{eq:modulo4list}
B_1, B_7, B_9, B_{14},B_{18}, B_{22}, B_{24}, B_{27}, B_{28}, B_{29}, B_{30}, B_{32}, B_{33}, B_{34}.
\end{gather*}
We note that this is the list of bases whose monodromy groups are contained in the alternating group. For other work observing congruence modulo four, see \cite{Modulo4,Modulo42}.

\subsection{All real or none real}
The bases
\begin{gather*}
B_{1},B_{3},B_{4},B_{7},B_{9},B_{18}, B_{24},B_{27},B_{28}
\end{gather*}
all have the property that if there are any real solutions, then all solutions are real. One reason for this could be that the discriminant of the branched cover does not intersect the subset of the parameter space we are sampling from, namely $[0,1]^6$.

\subsection{Real, positive, and realizable coincide}
For several of our experiments, we found that points in a fibre of $\pi_B$ over $[0,1]^6$ were real if and only if they were positive. This holds for the bases
\begin{gather*}
B_{1},B_{3},B_{7},B_8,B_9,B_{11},B_{12},B_{18},B_{19},B_{21},B_{22},
B_{23}.
\end{gather*}
Additionally, the following bases are those for which this coincidence extends to realizability
\begin{gather*}
B_{14},B_{15},B_{16},B_{17},B_{24},B_{29},B_{30},B_{31},B_{32},B_{33},B_{34},
B_{35}.
\end{gather*}
One can see this phenomenon explicitly with $B_{3}$: a fibre is completely determined by determining the only missing edge length, say $x_{23}$. Indeed, if one of the two possibilities for $x_{23}$ is real, the other must be as well since non-real solutions occur in conjugate pairs. We claim that $x_{23}$ is real if and only if it is positive. To see this, observe that 
\[
x_{23} = x_{12}+x_{13} \pm 2 \sqrt{x_{12}x_{13}-4x_{123}}
\]
as in Example \ref{ex:M2Basis} and so a real value of $x_{23}$ is negative if and only if 
$(x_{12}+x_{13})^2 < 4(x_{12}x_{13}-4x_{123})$ which is equivalent to $(x_{12}-x_{13})^2 < -16x_{123}$. This last inequality is never satisfied by positive $x_{12},x_{13}$ and $x_{123}$. Moreover, the Gram matrix of Theorem \ref{thm:Shoenberg} is 
\[
\begin{bmatrix}
2x_{12} & x_{12}+x_{13}-x_{23} & x_{12}+x_{14}-x_{24}\\
x_{12}+x_{13}-x_{23}  & 2x_{13} & x_{13}+x_{14}-x_{34} \\
 x_{12}+x_{14}-x_{24}& x_{13}+x_{14}-x_{34} & 2x_{14}
\end{bmatrix}.
\]
Notice that the upper left $2 \times 2$ minor is equal to $2x_{12}x_{13}+2x_{12}x_{23}+2x_{13}x_{23}-x_{12}^2-x_{13}^2-x_{23}^2$, or equivalently by Heron's formula, $x_{123}$. Hence, for any of our parameter values, this minor, and the other principal $2 \times 2$ minors, are always positive. The $3 \times 3$ determinant, on the other hand, need not always be positive for real values of $x_{12},\ldots,x_{34}$, as illustrated by the two vectors of squared edge lengths 
\begin{align*}
(x_{12},x_{13},x_{14},x_{23},x_{24},x_{34})  \in \{
 &(40.0, 90.0, 80.0, 16.862915010152395, 40.0, 30.0),\\ &
 (40.0, 90.0, 80.0, 243.1370849898476, 40.0, 30.0)\}.
\end{align*}
The second vector is not realizable by a $3$-simplex.

\subsection{Quantum Gravity, The Area-Length System, and Lorentzian $4$-simplices}
\label{secsec:physics}
General relativity uses a geometric description of spacetime via metric variables which describe lengths of curves in a spacetime manifold. Certain approaches to quantum gravity, on the other hand, uses area variables as their primary degrees of freedom as in loop quantum gravity \cite{Ashtekar:1986yd} and spin foam models \cite{Perez:2012wv}. Discretizing spacetime involves triangulating a $4$-dimensional Lorentzian manifold, and so it is of fundamental importance to understand the connection between the areas of $2$-faces of $\Delta_4$ and the edge lengths. 

\begin{figure}[!htpb]
\includegraphics[scale=0.4]{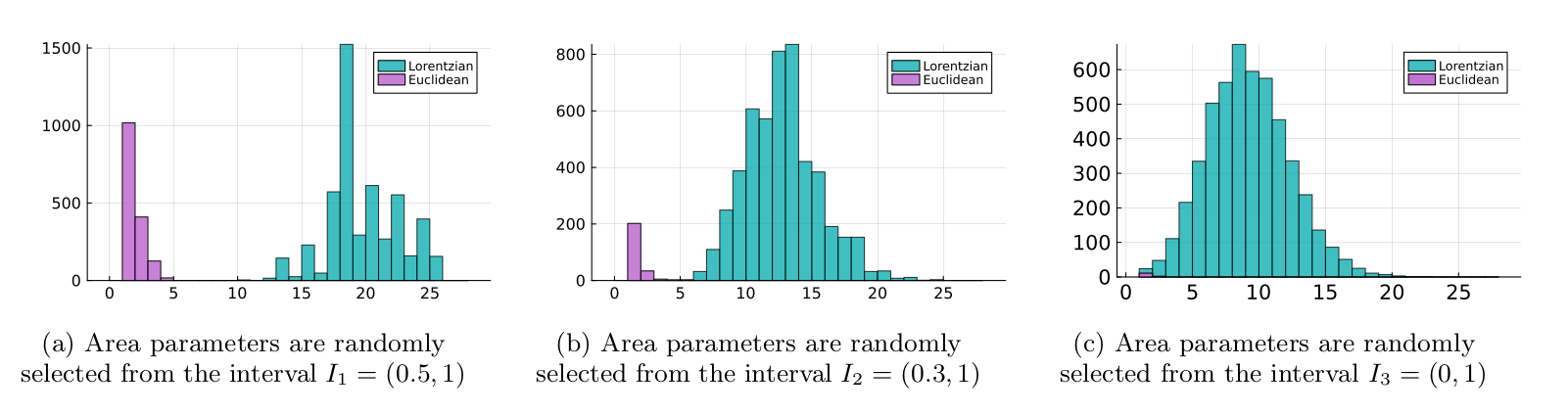}
\caption{Histograms indicating how many of the $64$ points in a fibre of $\pi_{B_{\textrm{area}}}$ are realizable as Euclidean or Lorentzian simplices.}
\label{fig:physics}
\end{figure}

This connection is encoded via the basis $B_{\textrm{area}}$ of $\mathcal M_4$ consisting of the $10$ triangular area variables. The corresponding branched cover has degree $64$ and Galois group $\ZZZ \wr \mathfrak S_{32}$ thanks to the involution $\textbf{x}_{ij} \mapsto -\textbf{x}_{ij}$ on the edge variables which fixes the corresponding polynomial system. This branched cover is the main topic of \cite{SethTaylor}. Included in that work are the following experiments, where we sample triangular areas uniformly from various intervals, and tabulate how many of the $64$ solutions are realizable as Euclidean or Lorentzian simplices. This data is summarized in Figure \ref{fig:physics}

\section*{Acknowledgements}
TB and MH are supported by an NSERC Discovery grant (RGPIN-2023-03551). SKA is supported by the Alexander von Humboldt foundation. The authors are grateful to Zvi Rosen and Daniel Bernstein for helpful conversations.

\newpage

{\small
\bibliographystyle{abbrv}
\bibliography{MyBibliography.bib}
}

\end{document}